\documentclass[10pt]{article}
\usepackage{amsmath, amsthm, amssymb, amsfonts, mathrsfs, mathtools}
\usepackage{graphicx}
\usepackage{makeidx}
\usepackage[left=2cm,right=2cm,top=2cm,bottom=2cm]{geometry}
\usepackage{caption}
\usepackage{subcaption}
\usepackage[section]{placeins}
\usepackage{enumitem}
\usepackage{tikz}
\usepackage{stmaryrd}

\usetikzlibrary{decorations.pathreplacing}
\tikzstyle{vertex}=[auto=left,circle,draw=black,fill=white, inner sep=1.5]

\usepackage{hyperref}
\hypersetup{colorlinks=true, citecolor={blue}, linkcolor=blue, filecolor=magenta, urlcolor=cyan}

\newtheorem{theorem}{Theorem}[section]

\newtheorem{lema}[theorem]{Lemma}
\newtheorem{corollary}{Corollary}[theorem]

\title{Laplacian Pair State Transfer on Edge Coronas}

\author{ Akash Kalita and Bikash Bhattacharjya\\
Department of Mathematics\\
Indian Institute of Technology Guwahati, India\\
akash.kalita@iitg.ac.in, b.bikash@iitg.ac.in}

\date{}
\begin{document}
\maketitle

\vspace{-0.3in}

\begin{center}{\textbf{Abstract}}\end{center}
Let $G$ be a graph with edge set $E(G) = \{e_1, \ldots, e_m\}$, and let $\mathcal{H} = (H_1, \ldots, H_m)$ be an $m$-tuple of graphs. The edge corona $G\diamond\mathcal{H}$ is obtained from the disjoint union of $G,H_1,\ldots,H_m$ by joining both end vertices of each edge $e_j$ to every vertex of $H_j$ for $1\leq j\leq m$. In this paper, we investigate Laplacian perfect pair state transfer and Laplacian pretty good pair state transfer on edge coronas. We establish necessary and sufficient conditions for the existence of Laplacian perfect pair state transfer. We also derive sufficient conditions for the non-existence of Laplacian perfect pair state transfer. Further, we 
identify conditions under which the absence of Laplacian perfect pair state transfer guarantees the existence of Laplacian pretty good pair state transfer. Additionally, we develop constructions based on Cartesian product and tensor product to obtain a few more infinite families of edge coronas admitting Laplacian pretty good pair state transfer. Finally, we investigate the existence of Laplacian pretty good pair state transfer on the edge corona of complete graphs.


\section{Introduction}
Throughout this paper, we only consider finite simple undirected connected graphs. Let $G$ be a graph on $n$ vertices and $M_G$ be a real symmetric matrix associated with $G$ such as the adjacency matrix or the Laplacian matrix or the signless Laplacian matrix. The \textit{transition matrix} of $G$ is denoted by $U_G(t)$ and is defined as 
\begin{align*}
U_G(t) :=  \exp(-\mathbf{i}t M_G) = \sum_{j = 0}^{\infty} \frac{(-\mathbf{i}t M_G)^j}{j!},
\end{align*}
where $\mathbf{i} := \sqrt{-1}$ and $t$ is a real variable. For a vertex $a$ of $G$, the vertex state associated with $a$ is defined as the standard basis vector $\textbf{e}_a$. The \textit{pair state} associated with a pair of vertex $(a_1,a_2)$ is defined as the vector $\frac{1}{\sqrt{2}}\left(\textbf{e}_{a_1} - \textbf{e}_{a_2}
\right)$. For convenience, we consider $\textbf{e}_{a_1} - \textbf{e}_{a_2}$ as the pair state associated with $(a_1,a_2)$.

Let $\textbf{x}, \textbf{y} \in {\mathbb{R}}^n$ be two linearly independent unit vectors. A graph $G$ is said to admit \textit{perfect real state transfer} from $\textbf{x}$ to $\textbf{y}$ if there exists a non-zero real number $t$ such that $\left|{\textbf{x}}^TU_G(t){\textbf{y}}\right| = 1$. If $\textbf{x} = \textbf{e}_{a_1}$ and $\textbf{y} = \textbf{e}_{a_2}$, where $a_1$ and $a_2$ are vertices of $G$, then perfect real state transfer is called perfect vertex state transfer. If $M_G=L_G$, $\textbf{x} = \textbf{e}_{a_1} - \textbf{e}_{a_2}$ and $\textbf{y} = \textbf{e}_{b_1} - \textbf{e}_{b_2}$, where $a_1,a_2,b_1,b_2$ are vertices of $G$, then perfect real state transfer is called Laplacian perfect pair state transfer (Pair-LPST).  

A graph $G$ is said to admit \textit{pretty good real state transfer} from $\textbf{x}$ to $\textbf{y}$ if for every positive real number $\varepsilon$, there exists a non-zero real number $t$ such that 
$\left|\left|{\textbf{x}}^TU_G(t){\textbf{y}}\right|-1\right|<\varepsilon$. If $\textbf{x} = \textbf{e}_{a_1}$ and $\textbf{y} = \textbf{e}_{a_2}$, where $a_1$ and $a_2$ are vertices of $G$, then pretty good real state transfer is called pretty good vertex state transfer. If $M_G=L_G$, $\textbf{x} = \textbf{e}_{a_1} - \textbf{e}_{a_2}$ and $\textbf{y} = \textbf{e}_{b_1} - \textbf{e}_{b_2}$, where $a_1,a_2,b_1,b_2$ are vertices of $G$, then pretty good real state transfer is called Laplacian pretty good pair state transfer (Pair-LPGST).  

Let $G$ be a graph with the vertex set $V(G) = \{v_1, \ldots, v_n\}$ and the edge set $E(G) = \{e_1, \ldots, e_m\}$. Also, let $\mathcal{H}$ be an $n$-tuple of graphs such that $\mathcal{H} = (H_1, \ldots, H_n)$. The \textit{vertex corona} of $G$ and $\mathcal{H}$ is denoted by $G \circ \mathcal{H}$ and is obtained by taking the disjoint union of $G, H_1,\dots,H_n$ and then joining the end vertices of $v_i$ to every vertex of $H_i$ for $1 \leq i \leq n$. In the definition of vertex corona, if $H_1 = \cdots = H_n = H$, then it is called \textit{homogeneous} vertex corona and is denoted by $G \circ H$. Let $\mathcal{H}$ be an $m$-tuple of graphs such that $\mathcal{H} = (H_1, \ldots, H_m)$. The \textit{edge corona} of $G$ and $\mathcal{H}$ is denoted by $G \diamond \mathcal{H}$ and is obtained by taking the disjoint union of $G, H_1,\ldots,H_m$ and then joining the end vertices of $e_j$ to every vertex of $H_j$ for $1 \leq j \leq m$. In the definition of edge corona, if $H_1 = \cdots = H_m = H$, then it is called \textit{homogeneous} edge corona and is denoted by $G \diamond H$.   

In~\cite{Chen2020}, Chen and Godsil introduced the notion of Pair-LPST on graphs. They 
provided a theoretical introduction to Pair-LPST and also provided methods to construct infinite families of graphs admitting Pair-LPST. Cao~\cite{Xiwang2021}, Cao and Wan~\cite{XiwangCao2022} and Luo et al.~\cite{GaojunLuo2022} studied Pair-LPST on cubelike graphs, abelian Cayley graphs and normal Cayley graphs over dihedral groups, respectively. Jiang et al.~\cite{PairStateTransferQGraph2025} investigated the existence of Pair-LPST and Pair-LPGST on $\mathcal{Q}$-graphs of regular graphs. Recently, Jiang et al.~\cite{DoubleCover2026} constructed infinite families of graphs admitting Pair-LPST using graph operation tensor product and double cover. 

In the last decade, several researchers studied vertex state transfer on different variants of corona. In the year 2016, Ackelsberg et al.~\cite{Ackelsberg2016} proved that the vertex corona $G \circ \mathcal{H}$ does not admit Laplacian perfect vertex state transfer if $G$ contains at least two vertices. In contrast, they prove that $G \circ \mathcal{H}$ can admits Laplacian pretty good vertex state transfer and provided the first infinite families of graphs admitting Laplacian pretty good vertex state transfer. Following their work, Ackelsberg et al.~\cite{Ackelsberg2017} and Tian et al.~\cite{Tian2021} investigated the existence of vertex state transfer on $G \circ \mathcal{H}$ with respect to the adjacency matrix and the signless Laplacian matrix, respectively. Li et al.~\cite{Edgecoronastatetransfer2022} studied the existence of Laplacian vertex state transfer on edge coronas. Wang and Liu~\cite{Wang2021}, and Wang and Liu~\cite{Wang2022} studied Laplacian vertex state transfer on edge complemented coronas and vertex complemented coronas, respectively. Compared to the vertex state case, there are very few articles that consider the existence of Pair-LPST and Pair-LPGST on coronas. In~\cite{VertexcoronaPairstatetransfer2023}, Wang et al. studied Pair-LPST and Pair-LPGST on homogeneous vertex coronas. Wang et al.~\cite{VertexcoronaPairstatetransfer2023} provided some sufficient conditions for the existence of Pair-LPGST on homogeneous vertex coronas. However, they did not construct infinite families admitting Pair-LPGST. Recently, Monterde and Pal~\cite{MonterdePal2026} obtain some results about Pair-LPST and Pair-LPGST on homogeneous vertex coronas. They also mentioned that analogous results can be obtained for other variants of coronas including the homogeneous edge corona. Despite all of this, Laplacian pair state transfer on other variants of corona (not necessarily homogeneous) remains largely unexplored. This inspired us to do a systematic study of  
Pair-LPST and Pair-LPGST on edge coronas (not necessarily homogeneous).

The paper is organized as follows. Section~\ref{preliminaries} is devoted to some basic definitions and a couple of preliminary results which are required to prove the main results of this paper. In Section~\ref{Pair-LPST on edge corona}, we obtain some necessary and sufficient conditions for the existence of Pair-LPST on edge coronas. Furthermore, we provide some sufficient conditions for the non-existence of Pair-LPST on edge coronas. In Section~\ref{Pair-LPGST on edge corona}, we provide some sufficient conditions for the existence of Pair-LPGST on edge coronas. In Section~\ref{Infinite families}, we obtain several infinite families of edge coronas admitting Pair-LPGST using the sufficient conditions obtain in Section~\ref{Pair-LPGST on edge corona}. Additionally, we show that one can obtain more infinite families of edge coronas admitting Pair-LPGST utilizing Cartesian product and tensor product of graphs. It is interesting that we obtain infinitely many graphs admitting Pair-LPGST that fails to admit Pair-LPST. Finally, in Section~\ref{Edge corona of the complete graph} we investigate the existence of Pair-LPGST on the edge corona $K_{n_1} \diamond K_{n_2}$, where $n_1$ and $n_2$ are positive integers such that $n_1 \geq 4$.

\section{Preliminaries}\label{preliminaries}
In order to prove our main results of this paper, we need a couple of preliminary results, which can be found in the literature. We provide these results in this section. Let $m$ and $n$ be two positive integers. In what follows, $I_{n \times n}$ denotes the identity matrix of order $n$, $\textbf{0}_{m \times n}$ denotes the $m \times n$ zero matrix,  
$J_{m \times n}$ denotes the $m \times n$ all one matrix, $\textbf{1}_{n}$ denotes the all one column matrix with $n$ entries and $\textbf{0}_{n}$ denotes the zero column matrix with $n$ entries. Further, $\mathrm{Cay}(\Gamma, S)$ denotes a Cayley graph over the group $\Gamma$ with connection set $S$. In particular, we have $K_n = \mathrm{Cay}(\mathbb{Z}_n, \mathbb{Z}_n \setminus \{0\})$ and $Q_n = \mathrm{Cay}(\mathbb{Z}^n_2, \{\textbf{e}_1,\ldots,\textbf{e}_n\})$. The \textit{dihedral group} is denoted by $D_{2n}$ and is defined as $D_{2n} := \langle a,b \colon a^{n} = b^2 = \textbf{1}, bab = a^{-1}\rangle$, the \textit{semi-dihedral group} is denoted by $SD_{8n}$ and is defined as $SD_{8n} := \langle a,b \colon a^{4n} = b^2 = \textbf{1}, bab = a^{2n-1}\rangle$, and the \textit{dicyclic group} is denoted by $T_{4n}$ and is defined as $T_{4n} := \langle a,b \colon a^{2n} = \textbf{1}, a^n = b^2, b^{-1}ab = a^{-1}\rangle$.

Let $G$ be a graph on $n$ vertices and $L_G$ be the Laplacian matrix of $G$. Throughout this paper, the spectrum of $G$ is the spectrum of $L_G$, denoted Spec($G$), unless stated otherwise. Let the spectral decomposition of $L_G$ be given by
$$L_G = \sum_{\theta \in \mathrm{Spec}(G)}\theta F_{\theta}{(G)}.$$
Let $(a_1,a_2)$ and $(b_1,b_2)$ be two distinct pairs of $G$. The \textit{eigenvalue support} of $(a_1,a_2)$ is denoted by $\Phi_{a_1a_2}$ and is defined as $\Phi_{a_1a_2}:=\{\theta \colon F_{\theta}{(G)}(\textbf{e}_{a_1} - \textbf{e}_{b_1}) \neq \textbf{0}_n \}$. The pairs $(a_1,a_2)$ and $(b_1,b_2)$ are called \textit{strongly cospectral} if $F_{\theta}{(G)}(\textbf{e}_{a_1} - \textbf{e}_{a_2}) = \pm F_{\theta}{(G)}(\textbf{e}_{b_1} - \textbf{e}_{b_2})$ for all $\theta \in \mathrm{Spec}(G)$. If $(a_1,a_2)$ and $(b_1,b_2)$ are strongly cospectral, then we consider the sets $\Phi^{+}_{a_1a_2,b_1b_2}$ and $\Phi^{-}_{a_1a_2,b_1b_2}$ which are given by 
\begin{align*}
\Phi^{+}_{a_1a_2,b_1b_2} &= \{\theta \colon F_{\theta}{(G)}(\textbf{e}_{a_1} - \textbf{e}_{a_2}) = F_{\theta}{(G)}(\textbf{e}_{b_1} - \textbf{e}_{b_2}) \neq \textbf{0}_n\}~\mathrm{and}~\\
\Phi^{-}_{a_1a_2,b_1b_2} &= \{\theta \colon F_{\theta}{(G)}(\textbf{e}_{a_1} - \textbf{e}_{a_2}) = -F_{\theta}{(G)}(\textbf{e}_{b_1} - \textbf{e}_{b_2}) \neq \textbf{0}_n\}.
\end{align*}
It is clear that if the pairs $(a_1,a_2)$ and $(b_1,b_2)$ are strongly cospectral, then we have  
$$\Phi_{a_1a_2} = \Phi_{b_1b_2} = \Phi^{+}_{a_1a_2,b_1b_2} \cup \Phi^{-}_{a_1a_2,b_1b_2}.$$ 

The pair $(a_1,a_2)$ is called \textit{fixed} if $\Phi_{a_1a_2}$ is a singleton set. The following lemma follows from Godsil et al.~\cite{RealStateTransfer2025}. 
\begin{lema}\emph{\cite{RealStateTransfer2025}}\label{fixed state}
Let the pairs $(a_1,a_2)$ and $(b_1,b_2)$ be strongly cospectral in a graph $G$. Then both of $\Phi^{+}_{a_1a_2,b_1b_2}$ and $\Phi^{-}_{a_1a_2,b_1b_2}$ are non-empty, that is, $(a_1,a_2)$ and $(b_1,b_2)$ cannot be fixed.
\end{lema}

In~\cite{Chen2020}, Chen and Godsil provided necessary and sufficient conditions for the existence of Pair-LPST which we include in the form of the following two theorems.

\begin{theorem}\emph{\cite{Chen2020}}\label{Conditions for the existence of Pair-LPST}
Let $(a_1,a_2)$ and $(b_1,b_2)$ be two distinct pairs of a graph $G$. Further, let $\Phi_{a_1a_2}$ be the set $\{\theta_0, \ldots, \theta_{i}\}$ such that $\theta_0 \in \Phi^{+}_{a_1a_2,b_1b_2}$. Then $G$ admits Pair-LPST from $(a_1,a_2)$ to $(b_1,b_2)$  if and only if the following three conditions hold.
\begin{enumerate}
\item[(i)] The pairs $(a_1,a_2)$ and $(b_1,b_2)$ are strongly cospectral.
\item[(ii)] Either $\Phi_{a_1a_2} \subseteq \mathbb{Z}$ or there is an integer $u$ and a square-free integer $\Delta$ strictly greater than one such that 
$$\Phi_{a_1a_2} \subseteq \left\{\frac{1}{2}(u + v \sqrt{\Delta}) \colon ~v~\mathrm{is~an~even~integer}\right \} ~~\mathrm{or}~~ \Phi_{a_1a_2} \subseteq \left\{\frac{1}{2}(u + v \sqrt{\Delta}) \colon ~v~\mathrm{is~an~odd~integer}\right \}.$$
\item[(iii)] Let \[g= \left\{ \begin{array}{rl} \gcd \left(\{\theta_0-\theta_j\}_{j=0}^{i}\right) & \textrm{ if } ~\Phi_{a_1a_2} \subseteq \mathbb{Z}\\
 \gcd \left(\left\{ \frac{\theta_0-\theta_j}{\sqrt{\Delta}}\right\}_{j=0}^{i}\right) & \textrm{ if }~\Phi_{a_1a_2} \subseteq \left\{\frac{1}{2}(u + v \sqrt{\Delta}) \colon ~v~\mathrm{is~an~integer}\right \},\\
\end{array}\right. \]
where $u$ is an integer and $\Delta$ is a square-free integer strictly greater than one. Then

$(a)~\theta_j \in \Phi^{+}_{a_1a_2,b_1b_2}$ if and only if 
\[ \left\{ \begin{array}{rl} \frac{\theta_0-\theta_j}{g}  \text{ is even}  & \textrm{ if } ~\Phi_{a_1a_2} \subseteq \mathbb{Z}\\
\frac{\theta_0-\theta_j}{g\sqrt{\Delta}} \text{ is even}  & \textrm{ if }~\Phi_{a_1a_2} \subseteq \left\{\frac{1}{2}(u + v \sqrt{\Delta}) \colon ~v~\mathrm{is~an~integer}\right \},\\
\end{array}\right. \]

$(b)~\theta_j \in \Phi^{-}_{a_1a_2,b_1b_2}$ if and only if 
\[ \left\{ \begin{array}{rl} \frac{\theta_0-\theta_j}{g}  \text{ is odd}  & \textrm{ if } ~\Phi_{a_1a_2} \subseteq \mathbb{Z}\\
\frac{\theta_0-\theta_j}{g\sqrt{\Delta}} \text{ is odd}  & \textrm{ if }~\Phi_{a_1a_2} \subseteq \left\{\frac{1}{2}(u + v \sqrt{\Delta}) \colon ~v~\mathrm{is~an~integer}\right \},\\
\end{array}\right.\]
\end{enumerate}
If the conditions $(i)$, $(ii)$ and $(iii)$ hold, then $G$ admits Pair-LPST at the minimum time $t_0$, where 
\[t_0 = \left\{ \begin{array}{cl} \frac{\pi}{g }  & \textrm{ if } ~\Phi_{a_1a_2} \subseteq \mathbb{Z}\\
 \frac{\pi}{g \sqrt{\Delta}} & \textrm{ if }~\Phi_{a_1a_2} \subseteq \left\{\frac{1}{2}(u + v \sqrt{\Delta}) \colon ~v~\mathrm{is~an~integer}\right \}.\\
\end{array}\right. \]
\end{theorem}

\begin{theorem}\emph{\cite{Chen2020}}\label{More conditions for the existence of perfect pair state transfer}
Let $(a_1,a_2)$ and $(b_1,b_2)$ be two distinct pairs of a graph $G$, and $\theta_0 \in \Phi^{+}_{a_1a_2,b_1b_2}$. Then $G$ admits Pair-LPST from $(a_1,a_2)$ to $(b_1,b_2)$ at time $t_1$ if and only if the following three conditions hold.
\begin{enumerate}
\item[(i)] The pairs $(a_1,a_2)$ and $(b_1,b_2)$ are strongly cospectral.
\item[(ii)] There exists $k_1 \in \mathbb{Z}$ such that $t_1(\theta_0 - \theta) = 2k_1 \pi$ for all $\theta \in \Phi^{+}_{a_1a_2,b_1b_2}$. 
\item[(iii)] There exists $k_2 \in \mathbb{Z}$ such that $t_1(\theta_0 - \theta) = (2k_2+1) \pi$ for all $\theta \in \Phi^{-}_{a_1a_2,b_1b_2}$. 
\end{enumerate}
\end{theorem}
The next corollary follows directly from Lemma~\ref{fixed state}, Theorem~\ref{Conditions for the existence of Pair-LPST} and Theorem~\ref{More conditions for the existence of perfect pair state transfer}.

\begin{corollary}\label{minimum time Corollary}
Let $G$ admit Pair-LPST from $(a_1,a_2)$ to $(b_1,b_2)$ at the minimum time $t_0$. Also, let $G$ admit Pair-LPST from $(a_1,a_2)$ to $(b_1,b_2)$ at time $t_1$. Then $t_1$ must be a rational multiple of $t_0$. 
\end{corollary}

The next result is known as \textit{Kronecker's approximation theorem} which is important in the investigation of Pair-LPGST on graphs.
\begin{theorem}\emph{\cite{Apostol1990}}\label{Kronecker approximation theorem}
Let $\beta_1, \ldots, \beta_n$ be arbitrary real numbers. Also, let $1, v_1, \ldots, v_n$ be real and linearly independent over $\mathbb{Q}$. Then for every positive real number $\varepsilon$, there exist integers $\ell,s_1, \ldots, s_n$ such that 
$$|\ell v_i - s_i -\beta_i |<\varepsilon~~\mathrm{for}~1 \leq i \leq n.$$
\end{theorem}

In our applications of Kronecker's approximation theorem, we will use the following lemma to identify sets of numbers which are linearly independent over $\mathbb{Q}$.
\begin{lema}\emph{\cite{Richard1974}}\label{Linearly independent corollary}
The set $\{\sqrt{\Delta} \colon \Delta~\mathrm{is~a~square~free~integer}\}$ is linearly independent over $\mathbb{Q}$.
\end{lema}

Let $G$ be an $r$-regular connected graph with $r \geq 2,~V(G)=\{v_1,\ldots,v_{n_1}\}$ and $E(G)=\{e_1,\ldots,e_m\}$. Also, let $\mathcal{H}=(H_1,\ldots,H_m)$, where $H_j$ is connected and $|V(H_j)|=n_2$ for $1 \leq j \leq m$. Then the vertex set $V(G \diamond \mathcal{H})$ of the graph $G \diamond \mathcal{H}$ can be labelled as follows:
$$V(G \diamond \mathcal{H}) = \{(v_i,0) \colon 1 \leq i \leq n_1\} \cup \{(w,e_j) \colon w \in V(H_j)~\mathrm{and}~e_j \in E(G)~\mathrm{for}~1 \leq j \leq m\}.$$
For a vertex $(v_i,0)$ of $G \diamond \mathcal{H}$, the vertex state of $(v_i,0)$ is  denoted by $\textbf{e}_{(v_i,0)}$ and is defined as
$$\textbf{e}_{(v_i,0)} = \begin{pmatrix}
\textbf{e}_{v_i} \\
    \textbf{0}_{mn_2}     
\end{pmatrix}.$$
For a vertex $(w,e_j)$ of $G \diamond \mathcal{H}$, the vertex state of $(w,e_j)$ is  denoted by $\textbf{e}_{(w,e_j)}$ and is defined as
$$\textbf{e}_{(w,e_j)} = \begin{pmatrix}
\textbf{0}_{n_1} \\
\textbf{e}_j \otimes \textbf{e}_w  
\end{pmatrix}.$$

Let $M_G$ be the $(0,1)$ incidence matrix of $G$ such that $\{\textbf{y}_1,\ldots,\textbf{y}_t\}$ is an orthonormal basis of the null-space of $M_G$. Then it is well known that $t = m - n_1$ if $G$ is non-bipartite and $t = m - n_1 + 1$ if $G$ is bipartite. The following theorem gives the eigenvalues and the eigenprojectors of $G \diamond \mathcal{H}$, where $G$ is non-bipartite.
\begin{theorem}\emph{\cite{Edgecoronastatetransfer2022}}\label{Non-bipartite Spectrum}  
If $G$ is non-bipartite, the eigenvalues and the corresponding eigenprojectors of $G \diamond \mathcal{H}$ are as follows:
\begin{enumerate}
\item[(i)] Let $||\textbf{y}_j||$ denote the standard norm of $\textbf{y}_j$ in $\mathbb{R}^m$ for $1 \leq j \leq m-n_1$. Then 2 is an eigenvalue of $G \diamond \mathcal{H}$ with the corresponding eigenprojector given by
$$F_2(G \diamond \mathcal{H})=\begin{pmatrix}
         \textbf{0}_{n_1 \times n_1} &   \textbf{0}_{n_1 \times mn_2}\\
           \textbf{0}_{mn_2 \times n_1} & \sum_{j=1}^{m-n_1}\frac{1}{||\textbf{y}_j||^2}(\textbf{y}_j{\textbf{y}_j}^T) \otimes \frac{1}{n_2}J_{n_2 \times n_2}
       \end{pmatrix}.$$
\item[(ii)] Let $\mu$ be a non-zero eigenvalue of $H_j$ with the corresponding eigenprojector $F_{\mu}(H_j)$ for $1 \leq j \leq m$. Then $\mu+2$ is an eigenvalue of $G \diamond \mathcal{H}$ with the corresponding eigenprojector given by
$$F_{\mu+2}(G \diamond \mathcal{H})=\begin{pmatrix}
         \textbf{0}_{n_1 \times n_1} &   \textbf{0}_{n_1 \times mn_2}\\
           \textbf{0}_{mn_2 \times n_1} & \sum_{j=1}^{m}(\textbf{e}_j{\textbf{e}_j}^T) \otimes F_{\mu}(H_j)
       \end{pmatrix}.$$
\item[(iii)] Let $\theta$ be an eigenvalue of $G$ with the corresponding eigenprojector $F_{\theta}(G)$. Also, let 
$$\theta^{\pm}=\frac{rn_2+\theta+2\pm\sqrt{(rn_2+\theta-2)^2+4n_2(2r-\theta)}}{2}.$$
Then $\theta^{\pm}$ are eigenvalues of $G \diamond \mathcal{H}$ with the corresponding eigenprojector given by 
$$F_{\theta^{\pm}}(G \diamond \mathcal{H}) = \frac{(2-\theta^{\pm})^2}{(2-\theta^{\pm})^2 + (2r-\theta)n_2}\begin{pmatrix}
F_{\theta}(G) & \frac{1}{2-\theta^{\pm}}(F_{\theta}(G)M_G) \otimes \textbf{1}^T_{n_2}\\
\frac{1}{2-\theta^{\pm}}({M_G}^TF_{\theta}(G)) \otimes \textbf{1}_{n_2} & \frac{1}{(2-\theta^{\pm})^2}({M_G}^TF_{\theta}(G)M_G) \otimes J_{n_2 \times n_2}
\end{pmatrix}.$$
\end{enumerate}
\end{theorem}

The next theorem gives the eigenvalues and the eigenprojectors of $G \diamond \mathcal{H}$, where $G$ is bipartite.
\begin{theorem}\emph{\cite{Edgecoronastatetransfer2022}}\label{Bipartite Spectrum} 
If $G$ is bipartite with the bipartition $V_1 \cup V_2$, the eigenvalues and the corresponding eigenprojectors of $G \diamond \mathcal{H}$ are as follows:
\begin{enumerate}
\item[(i)] Let $||\textbf{y}_j||$ denote the standard norm of $\textbf{y}_j$ in $\mathbb{R}^m$ for $1 \leq j \leq m-n_1+1$. Then 2 is an eigenvalue of $G \diamond \mathcal{H}$ with the corresponding eigenprojector given by
$$F_2(G \diamond \mathcal{H})=\begin{pmatrix}
         \textbf{0}_{n_1 \times n_1} &   \textbf{0}_{n_1 \times mn_2}\\
           \textbf{0}_{mn_2 \times n_1} & \sum_{j=1}^{m-n_1+1}\frac{1}{||\textbf{y}_j||^2}(\textbf{y}_j{\textbf{y}_j}^T) \otimes \frac{1}{n_2}J_{n_2 \times n_2}
       \end{pmatrix}.$$
\item[(ii)] Let $\mu$ be a non-zero eigenvalue of $H_j$ with the corresponding eigenprojector $F_{\mu}(H_j)$ for $1 \leq j \leq m$. Then $\mu+2$ is a eigenvalue of $G \diamond \mathcal{H}$ with the corresponding eigenprojector given by
$$F_{\mu+2}(G \diamond \mathcal{H})=\begin{pmatrix}
         \textbf{0}_{n_1 \times n_1} & \textbf{0}_{n_1 \times mn_2}\\
          \textbf{0}_{mn_2 \times n_1} & \sum_{j=1}^{m}(\textbf{e}_j{\textbf{e}_j}^T) \otimes F_{\mu}(H_j)
       \end{pmatrix}.$$
\item[(iii)] Let $\theta$ be an eigenvalue of $G$ with the corresponding eigenprojector $F_{\theta}(G)$ such that $\theta \neq 2r$. Also, let 
$$\theta^{\pm}=\frac{rn_2+\theta+2\pm\sqrt{(rn_2+\theta-2)^2+4n_2(2r-\theta)}}{2}.$$
Then $\theta^{\pm}$ are eigenvalues of $G \diamond \mathcal{H}$ with the corresponding eigenprojector given by 
$$F_{\theta^{\pm}}(G \diamond \mathcal{H}) = \frac{(2-\theta^{\pm})^2}{(2-\theta^{\pm})^2 + (2r-\theta)n_2}\begin{pmatrix}
F_{\theta}(G) & \frac{1}{2-\theta^{\pm}}(F_{\theta}(G)M_G) \otimes \textbf{1}^T_{n_2}\\
\frac{1}{2-\theta^{\pm}}({M_G}^TF_{\theta}(G)) \otimes \textbf{1}_{n_2} & \frac{1}{(2-\theta^{\pm})^2}({M_G}^TF_{\theta}(G)M_G) \otimes J_{n_2 \times n_2}
\end{pmatrix}.$$
\item[(iv)] Let $R$ denote the matrix 
$$\frac{1}{n_1}\begin{pmatrix}
 J_{|V_1| \times |V_1|} & -J_{|V_1| \times |V_2|}\\
-J_{|V_2| \times |V_1|} &  J_{|V_2| \times |V_2|}
\end{pmatrix}.$$
Then $r(n_2+2)$ is an eigenvalue of $G \diamond \mathcal{H}$ with the corresponding eigenprojector given by
$$F_{r(n_2+2)}(G \diamond \mathcal{H}) = \begin{pmatrix}
 R & \textbf{0}_{n_1 \times mn_2}\\
 \textbf{0}_{mn_2 \times n_1} & \textbf{0}_{mn_2 \times mn_2}
\end{pmatrix}.$$
\end{enumerate}
\end{theorem}

\section{Pair-LPST on edge coronas}\label{Pair-LPST on edge corona}
In this section, we investigate the existence of Pair-LPST on edge coronas. In what follows, $G$ denotes an $r$-regular connected graph such that $r \geq 2,~|V(G)|=n_1$ and $|E(G)|=m$. Further, $\mathcal{H}=(H_1,\ldots,H_m)$, where $H_j$ is connected and $|V(H_j)|=n_2$ for $1 \leq j \leq m$. If $H_1 = \ldots = H_m = H$, then the following theorem follows from Monterde and Pal~\cite{MonterdePal2026} [see Theorem 7].
\begin{theorem}\emph{\cite{MonterdePal2026}}\label{Graphs with clusters Pair-LPST} 
Let $e_s$ be a fixed edge of $G$ and $h_i$ be a vertex of $H$, where $1 \leq i \leq 4$. Then the graph $H$ admits Pair-LPST from $(h_1,h_2)$ to $(h_3,h_4)$
if and only if the edge corona $G \diamond H$ admits Pair-LPST from $((h_1,e_s),(h_2,e_s))$ to $((h_3,e_s),(h_4,e_s))$.
\end{theorem}
In the following, we extend Theorem~\ref{Graphs with clusters Pair-LPST} to the more general form of edge corona $G \diamond (H_1,\ldots,H_m)$, where $H_1,\ldots,H_m$  can be distinct. 
\begin{theorem}
Let $e_s$ be a fixed edge of $G$ and $h_i$ be a vertex of $H_s$, where $1 \leq i \leq 4$. Then the graph $H_s$ admits Pair-LPST from $(h_1,h_2)$ to $(h_3,h_4)$
if and only if the edge corona $G \diamond \mathcal{H}$ admits Pair-LPST from $((h_1,e_s),(h_2,e_s))$ to $((h_3,e_s),(h_4,e_s))$.
\end{theorem}
\begin{proof}
Let $G$ be non-bipartite. Then using Theorem~\ref{Non-bipartite Spectrum}, we obtain that
\begin{align}
F_2(G \diamond \mathcal{H})(\textbf{e}_{(h_1,e_s)} - \textbf{e}_{(h_2,e_s)}) = \textbf{0}_{n_1+mn_2} = F_2(G \diamond \mathcal{H})(\textbf{e}_{(h_3,e_s)} - \textbf{e}_{(h_4,e_s)})\label{F_1}~\mathrm{and}~\\ 
F_{\theta^{\pm}}(G \diamond \mathcal{H})(\textbf{e}_{(h_1,e_s)} - \textbf{e}_{(h_2,e_s)}) = \textbf{0}_{n_1+mn_2} = F_{\theta^{\pm}}(G \diamond \mathcal{H})(\textbf{e}_{(h_3,e_s)} - \textbf{e}_{(h_4,e_s)})\label{F_2},
\end{align}
for all $\theta \in \mathrm{Spec}(G)$. Further, for each non-zero eigenvalue $\mu$ of $H_s$, we obtain that
\begin{align}
&F_{\mu+2}(G \diamond \mathcal{H})(\textbf{e}_{(h_1,e_s)} - \textbf{e}_{(h_2,e_s)}) = \begin{pmatrix}
\textbf{0}_{n_1}\\
\textbf{e}_s \otimes F_{\mu}(H_s)(\textbf{e}_{h_1}-\textbf{e}_{h_2})
\end{pmatrix}\label{f_1}~\mathrm{and}~\\
&F_{\mu+2}(G \diamond \mathcal{H})(\textbf{e}_{(h_3,e_s)} - \textbf{e}_{(h_4,e_s)}) = \begin{pmatrix}
\textbf{0}_{n_1}\\
\textbf{e}_s \otimes F_{\mu}(H_s)(\textbf{e}_{h_3}-\textbf{e}_{h_4})
\end{pmatrix}\label{f_2}.
\end{align}
From Equations~(\ref{F_1}), (\ref{F_2}), (\ref{f_1}) and (\ref{f_2}), it is clear that 
\begin{align}
&\Phi_{(h_1,e_s)(h_2,e_s)}=\{\mu+2 \colon \mu \in \Phi_{h_1h_2}\}\label{A_1}\\
&\Phi^{+}_{(h_1,e_s)(h_2,e_s),(h_3,e_s)(h_4,e_s)}=\{\mu+2 \colon \mu \in \Phi^{+}_{h_1h_2,h_3h_4}\}\label{A_2}~\mathrm{and}~\\
&\Phi^{-}_{(h_1,e_s)(h_2,e_s),(h_3,e_s)(h_4,e_s)}=\{\mu+2 \colon \mu \in \Phi^{-}_{h_1h_2,h_3h_4}\}\label{A_3}.
\end{align}
Let $H_s$ admit Pair-LPST from $(h_1,h_2)$ to $(h_3,h_4)$. Then $(h_1,h_2)$ and $(h_3,h_4)$ are strongly cospectral. Now Equation~(\ref{f_1}) and Equation~(\ref{f_2}) together implies that $((h_1,e_s),(h_2,e_s))$ and $((h_3,e_s),(h_4,e_s))$ are strongly cospectral. Applying Equations~(\ref{A_1}), (\ref{A_2}) and (\ref{A_3}), the conditions $(ii)$ and $(iii)$ of Theorem~\ref{Conditions for the existence of Pair-LPST} can be easily verified. Therefore, $G \diamond \mathcal{H}$ admits Pair-LPST from $((h_1,e_s),(h_2,e_s))$ to $((h_3,e_s),(h_4,e_s))$.

Conversely, let the edge corona $G \diamond \mathcal{H}$ admit Pair-LPST from $((h_1,e_s),(h_2,e_s))$ to $((h_3,e_s),(h_4,e_s))$. Then $((h_1,e_s),(h_2,e_s))$ and $((h_3,e_s),(h_4,e_s))$ are strongly cospectral. Now Equation~(\ref{f_1}) and Equation~(\ref{f_2}) together implies that $(h_1,h_2)$ and $(h_3,h_4)$ are strongly cospectral. Also, the conditions $(ii)$ and $(iii)$ of Theorem~\ref{Conditions for the existence of Pair-LPST} can be easily verified. Therefore, $H_s$ admits Pair-LPST from $(h_1,h_2)$ to $(h_3,h_4)$. 

If $G$ is bipartite, then the proof is similar to that of the non-bipartite case, and hence the details are omitted to avoid repetitive argument. 
\end{proof}
Theorem~\ref{Graphs with clusters Pair-LPST} fails to give any information about the existence of Pair-LPST from the remaining pair (for instance from pair $((a_1,0), (a_2,0))$ to $((b_1,0), (b_2,0))$, where $a_1,a_2,b_1,b_2 \in V(G)$). Now we explore the existence of Pair-LPST on $G \diamond \mathcal{H}$ from the remaining pair. We first prove the following two lemmas borrowing some ideas from Jiang et al.~\cite{PairStateTransferQGraph2025} [see Theorem 3.1].
\begin{lema}\label{Laplacian perfect pair state transfer lemma non-bipartite}
Let $G$ be non-bipartite, $(a_1,0)$ and $(a_2,0)$ be two vertices of the edge corona $G \diamond \mathcal{H}$. Then 
\begin{enumerate}
\item[(a)] $\theta^{+} \in \Phi_{(a_1,0)(a_2,0)}$ if and only if $\theta^{-} \in \Phi_{(a_1,0)(a_2,0)}$ for each $\theta \in \mathrm{Spec}(G)$.
\item[(b)] If a pair is strongly cospectral with $((a_1,0),(a_2,0))$, then it must be of the form $((b_1,0),(b_2,0))$. Let $((a_1,0),(a_2,0))$ and $((b_1,0),(b_2,0))$ be strongly cospectral. Then $\theta^{+} \in \Phi^{+}_{(a_1,0)(a_2,0),(b_1,0)(b_2,0)}$ if and only if $\theta^{-} \in \Phi^{+}_{(a_1,0)(a_2,0),(b_1,0)(b_2,0)}$ for each $\theta \in \mathrm{Spec}(G)$. Further, $\theta^{+} \in \Phi^{-}_{(a_1,0)(a_2,0),(b_1,0)(b_2,0)}$ if and only if $\theta^{-} \in \Phi^{-}_{(a_1,0)(a_2,0),(b_1,0)(b_2,0)}$ for each $\theta \in \mathrm{Spec}(G)$.
\end{enumerate} 
\end{lema}
\begin{proof}
$(a)$ Let $2-\theta^{\pm } = 0$ for some $\theta \in \mathrm{Spec}(G)$. Since 
$$\theta^{\pm}=\frac{rn_2+\theta+2\pm\sqrt{(rn_2+\theta-2)^2+4n_2(2r-\theta)}}{2},$$
now it follows that $\theta = 2r$. This contradicts the fact that $G$ is non-bipartite. Thus $2-\theta^{\pm } \neq 0$ for each $\theta \in \mathrm{Spec}(G)$. Applying Theorem~\ref{Non-bipartite Spectrum}, for each $\theta \in \mathrm{Spec}(G)$ we obtain that 
\begin{equation}\label{x_1}
F_{\theta^{\pm }}(G \diamond \mathcal{H})(\textbf{e}_{(a_1,0)} - \textbf{e}_{(a_2,0)}) = \frac{{(2-{\theta}^{\pm})}^2}{(2-\theta^{\pm })^2+(2r-\theta)n_2}\begin{pmatrix}
        F_{\theta}(G)(\textbf{e}_{a_1} - \textbf{e}_{a_2})  \\
        \frac{1}{2-\theta^{\pm}}{M^T_G}F_{\theta}(G)(\textbf{e}_{a_1} - \textbf{e}_{a_2}) \otimes \textbf{1}_{n_2} 
       \end{pmatrix}.
\end{equation}  
From Equation~(\ref{x_1}), we find that $F_{\theta^{+}}(G \diamond \mathcal{H})(\textbf{e}_{(a_1,0)} - \textbf{e}_{(a_2,0)}) \neq \textbf{0}_{n_1+mn_2}$ if and only if $F_{\theta^{-}}(G \diamond \mathcal{H})(\textbf{e}_{(a_1,0)} - \textbf{e}_{(a_2,0)}) \neq \textbf{0}_{n_1+mn_2}$. Therefore $\theta^{+} \in \Phi_{(a_1,0)(a_2,0)}$ if and only if $\theta^{-}\in  \Phi_{(a_1,0)(a_2,0)}$ for each $\theta \in \mathrm{Spec}(G)$.

$(b)$ First we prove that $((a_1,0),(a_2,0))$ and $((b_1,0),(b_2,e_t))$ can never be strongly cospectral. If possible, let $((a_1,0),(a_2,0))$ and $((b_1,0),(b_2,e_t))$ be strongly cospectral. Then we have 
 \begin{equation}\label{y_2}
F_{\theta^{\pm }}(G \diamond \mathcal{H})(\textbf{e}_{(a_1,0)} - \textbf{e}_{(a_2,0)}) = \pm F_{\theta^{\pm }}(G \diamond \mathcal{H})(\textbf{e}_{(b_1,0)} - \textbf{e}_{(b_2,e_t)})~~~\mathrm{for~each}~\theta \in \mathrm{Spec}(G).
\end{equation}
Since $G$ is a connected graph on $n_1$ vertices, we have $F_{0}(G) = \frac{1}{n_1}J_{n_1 \times n_1}$. This implies $F_{0}(G)(\textbf{e}_{a_1} - \textbf{e}_{a_2})=\textbf{0}_{n_1}$. Also, we have $0^{+}=rn_2+2$. Now Equation~(\ref{x_1}) and Equation~(\ref{y_2}) yield that 
\begin{equation}\label{x_2}
F_{rn_2+2}(G \diamond \mathcal{H})(\textbf{e}_{(a_1,0)} - \textbf{e}_{(a_2,0)}) = \textbf{0}_{n_1+mn_2} = F_{rn_2+2}(G \diamond \mathcal{H})(\textbf{e}_{(b_1,0)} - \textbf{e}_{(b_2,e_t)}).
\end{equation}
Using Theorem~\ref{Non-bipartite Spectrum}, we obtain 
\begin{equation}\label{x_3}
F_{rn_2+2}(G \diamond \mathcal{H}) = \frac{rn_2}{rn_2+2}\begin{pmatrix}
         \frac{1}{n_1} J_{n_1\times n_1} &  \frac{-2}{rn_1n_2} J_{n_1 \times mn_2}\\
          \frac{-2}{rn_1n_2}J_{mn_2 \times n_1} &  \frac{4}{r^2n_1n^2_2} J_{mn_2 \times mn_2}
       \end{pmatrix}.
\end{equation}
From Equation~(\ref{x_3}), it follows that $F_{rn_2+2}(G \diamond \mathcal{H})(\textbf{e}_{(b_1,0)} - \textbf{e}_{(b_2,e_t)}) \neq \textbf{0}_{n_1+mn_2}$, a contradiction to Equation~(\ref{x_2}). Therefore, $((a_1,0),(a_2,0))$ and $((b_1,0),(b_2,e_t))$ can never be strongly cospectral.

From Equation~(\ref{x_3}), it is clear that if a pair is strongly cospectral with $((a_1,0),(a_2,0))$, then it is either of the form $((b_1,0),(b_2,0))$ or of the form $((b_1,e_s),(b_2,e_t))$. If possible, let $((a_1,0),(a_2,0))$ and $((b_1,e_s),(b_2,e_t))$ be strongly cospectral. Now
\begin{equation}\label{x_4}
F_{\theta^{\pm }}(G \diamond \mathcal{H})(\textbf{e}_{(b_1,e_s)} - \textbf{e}_{(b_2,e_t)}) = \frac{{(2-{\theta}^{\pm})}^2}{(2-\theta^{\pm })^2+(2r-\theta)n_2}\begin{pmatrix}
         \frac{1}{{(2-\theta^{\pm})}^2}F_{\theta}(G)M_G(\textbf{e}_s - \textbf{e}_t)  \\
         \frac{1}{2-\theta^{\pm}}{M^T_G}F_{\theta}(G)M_G (\textbf{e}_s - \textbf{e}_t) \otimes \textbf{1}_{n_2}
       \end{pmatrix}.
\end{equation}
Suppose that $\theta^{\pm} \in \Phi_{(a_1,0)(a_2,0)}$ for some $\theta \in \mathrm{Spec}(G)$. The following two cases arise. 

\noindent \textbf{Case 1.} Let $\theta^{+} \in \Phi^{+}_{(a_1,0)(a_2,0),(b_1,e_s)(b_2,e_t)}$ and $\theta^{-} \in \Phi^{+}_{(a_1,0)(a_2,0),(b_1,e_s)(b_2,e_t)}$. Then we have  
\begin{align}
&F_{\theta^{+ }}(G \diamond \mathcal{H})(\textbf{e}_{(a_1,0)} - \textbf{e}_{(a_2,0)}) = F_{\theta^{+}}(G \diamond \mathcal{H})(\textbf{e}_{(b_1,e_s)} - \textbf{e}_{(b_2,e_t)})~\mathrm{and}~\label{T_1}\\
&F_{\theta^{-}}(G \diamond \mathcal{H})(\textbf{e}_{(a_1,0)} - \textbf{e}_{(a_2,0)}) = -F_{\theta^{-}}(G \diamond \mathcal{H})(\textbf{e}_{(b_1,e_s)} - \textbf{e}_{(b_2,e_t)})\label{T_2}.
\end{align}
Using equations~(\ref{x_1}) and~(\ref{x_4}) in equations~(\ref{T_1}) and~(\ref{T_2}), we obtain  
\begin{align*}
&F_{\theta}(G)M_G(\textbf{e}_s - \textbf{e}_t) = (2-\theta^{+})F_{\theta}(G)(\textbf{e}_{a_1} - \textbf{e}_{a_2})~\mathrm{and}~\\
&F_{\theta}(G)M_G(\textbf{e}_s - \textbf{e}_t) = -(2-\theta^{-})F_{\theta}(G)(\textbf{e}_{a_1} - \textbf{e}_{a_2}).
\end{align*}
As $F_{\theta}(G)(\textbf{e}_{a_1} - \textbf{e}_{a_2}) \neq \textbf{0}_{n_1}$, we have from the preceding two equations that $\theta^{+} + \theta^{-} = 4$. This implies $rn_2+\theta=2$, this is a contradiction as $r \geq 2,~n_2$ is a positive integer and $\theta$ is positive. We get the same contradiction if  $\theta^{-} \in \Phi^{+}_{(a_1,0)(a_2,0),(b_1,e_s)(b_2,e_t)}$ and $\theta^{+} \in \Phi^{-}_{(a_1,0)(a_2,0),(b_1,e_s)(b_2,e_t)}$.

\noindent \textbf{Case 2.} Let $\theta^{+},\theta^{-} \in \Phi^{+}_{(a_1,0)(a_2,0),(b_1,e_s)(b_2,e_t)}$. In this case also, proceeding as in Case 1 we find $\theta^{+} = \theta^{-}$. This implies $(rn_2+\theta-2)^2+4n_2(2r-\theta)=0$. This is a contradiction as $0 < (rn_2+\theta-2)^2+4n_2(2r-\theta)$ for each $\theta \in \mathrm{Spec}(G)$. We get the same contradiction if $\theta^{+}, \theta^{-} \in \Phi^{-}_{(a_1,0)(a_2,0),(b_1,e_s)(b_2,e_t)}$.

Thus $\theta^{\pm} \notin \Phi_{(a_1,0)(a_2,0)}$ for each $\theta \in \mathrm{Spec}(G)$. By Theorem~\ref{Non-bipartite Spectrum}, the eigenvalues of $G \diamond \mathcal{H}$ are $2, \mu + 2$ and $\theta^{\pm}$. It is clear that $2, \mu + 2 \notin \Phi_{(a_1,0)(a_2,0)}$. However, By Lemma~\ref{fixed state}, $\left|\Phi_{(a_1,0)(a_2,0)}\right| \geq 2$. So, $\Phi_{(a_1,0)(a_2,0)}$ must contain both of $\theta^{+}$ and $\theta^{-}$ for some $\theta \in \mathrm{Spec}(G)$, which is not possible. Therefore, $((a_1,0),(a_2,0))$ and $((b_1,e_s),(b_2,e_t))$ can never be strongly cospectral.

Let $((a_1,0),(a_2,0))$ and $((b_1,0),(b_2,0))$ be strongly cospectral. Then we have
\begin{equation}\label{x_5}
F_{\theta^{\pm }}(G \diamond \mathcal{H})(\textbf{e}_{(b_1,0)} - \textbf{e}_{(b_2,0)}) = \frac{{(2-{\theta}^{\pm})}^2}{(2-\theta^{\pm })^2+(2r-\theta)n_2}\begin{pmatrix}
        F_{\theta}(G)(\textbf{e}_{b_1} - \textbf{e}_{b_2})  \\
        \frac{1}{2-\theta^{\pm}}{M^T_G}F_{\theta}(G)(\textbf{e}_{b_1} - \textbf{e}_{b_2}) \otimes \textbf{1}_{n_2} 
       \end{pmatrix}.
\end{equation}  
From equations (\ref{x_1}) and (\ref{x_5}), we find  $\theta^{+} \in \Phi^{+}_{(a_1,0)(a_2,0),(b_1,0)(b_2,0)}$ if and only if $F_{\theta}(G)(\textbf{e}_{a_1} - \textbf{e}_{a_2}) = F_{\theta}(G)(\textbf{e}_{b_1} - \textbf{e}_{b_2}) \neq \textbf{0}_{n_1}$. Similarly, $\theta^{-} \in \Phi^{+}_{(a_1,0)(a_2,0),(b_1,0)(b_2,0)}$ if and only if $F_{\theta}(G)(\textbf{e}_{a_1} - \textbf{e}_{a_2}) = F_{\theta}(G)(\textbf{e}_{b_1} - \textbf{e}_{b_2}) \neq \textbf{0}_{n_1}$. Hence, $\theta^{+} \in \Phi^{+}_{(a_1,0)(a_2,0),(b_1,0)(b_2,0)}$ if and only if $\theta^{-} \in \Phi^{+}_{(a_1,0)(a_2,0),(b_1,0)(b_2,0)}$ for each $\theta \in \mathrm{Spec}(G)$. In a similar way, $\theta^{+} \in \Phi^{-}_{(a_1,0)(a_2,0),(b_1,0)(b_2,0)}$ if and only if $\theta^{-} \in \Phi^{-}_{(a_1,0)(a_2,0),(b_1,0)(b_2,0)}$ for each $\theta \in \mathrm{Spec}(G)$. This completes the proof.    
\end{proof}

Following the same procedure as in the Lemma~\ref{Laplacian perfect pair state transfer lemma non-bipartite}, we prove the next lemma. The details are omitted.
\begin{lema}\label{Laplacian perfect pair state transfer lemma bipartite}
Let $G$ be bipartite, $(a_1,0)$ and $(a_2,0)$ be two vertices of the edge corona $G \diamond \mathcal{H}$. Then 
\begin{enumerate}
\item[(i)] $\theta^{+} \in \Phi_{(a_1,0)(a_2,0)}$ if and only if $\theta^{-} \in \Phi_{(a_1,0)(a_2,0)}$ for each $\theta \in \mathrm{Spec}(G) \setminus \{2r\}$.
\item[(ii)] If a pair is strongly cospectral with $((a_1,0),(a_2,0))$, then it must be of the from $((b_1,0),(b_2,0))$. Let $((a_1,0),(a_2,0))$ and $((b_1,0),(b_2,0))$ be strongly cospectral. Then $\theta^{+} \in \Phi^{+}_{(a_1,0)(a_2,0),(b_1,0)(b_2,0)}$ if and only if $\theta^{-} \in \Phi^{+}_{(a_1,0)(a_2,0),(b_1,0)(b_2,0)}$ for each $\theta \in \mathrm{Spec}(G) \setminus \{2r\}$. Further, $\theta^{+} \in \Phi^{-}_{(a_1,0)(a_2,0),(b_1,0)(b_2,0)}$ if and only if $\theta^{-} \in \Phi^{-}_{(a_1,0)(a_2,0),(b_1,0)(b_2,0)}$ for each $\theta \in \mathrm{Spec}(G) \setminus \{2r\}$.
\end{enumerate} 
\end{lema} 
Now we use Lemma~\ref{Laplacian perfect pair state transfer lemma non-bipartite} and Lemma~\ref{Laplacian perfect pair state transfer lemma bipartite} to prove the next theorem. Here we have the assumption that $G$ is an $r$-regular graph such that $r\geq3$.

\begin{theorem}\label{non-existence of Pair-LPST First Theorem}
Let $a_1$ and $a_2$ be two vertices of $G$. Also, let $(r-1)/(1+n_2)$ not be an integer. Then the pair $((a_1,0),(a_2,0))$ cannot be involved in Pair-LPST.  
\end{theorem}
\begin{proof}
First consider that $G$ is non-bipartite. If possible, let $((a_1,0),(a_2,0))$ be involve in Pair-LPST. Then from the proof of Lemma~\ref{Laplacian perfect pair state transfer lemma non-bipartite}, we obtain that $G \diamond \mathcal{H}$ admits Pair-LPST from $((a_1,0),(a_2,0))$ to $((b_1,0),(b_2,0))$, for some vertices $b_1$ and $b_2$ of $G$. From Theorem~\ref{Non-bipartite Spectrum}, it is clear that $\mu+2 \notin \Phi_{(a_1,0)(a_2,0)}$, where $\mu \in \mathrm{Spec}(H_j)$ for $1 \leq j \leq m$. Further, Equation~(\ref{x_2}) yields that $0^{\pm} \notin \Phi_{(a_1,0)(a_2,0)}$.
Thus we have $\Phi_{(a_1,0)(a_2,0)} \subseteq \{\theta^{\pm } \colon \theta \in \mathrm{Spec}(G)\setminus \{0\}\}$. By Theorem~\ref{Conditions for the existence of Pair-LPST}, either all the elements in $\Phi_{(a_1,0)(a_2,0)}$ are integers or all the elements in $\Phi_{(a_1,0)(a_2,0)}$ are quadratic integers. From Lemma~\ref{Laplacian perfect pair state transfer lemma non-bipartite}, we have either $\theta^{\pm } \in \Phi_{(a_1,0)(a_2,0)}$ or $\theta^{\pm } \notin \Phi_{(a_1,0)(a_2,0)}$ for each $\theta \in \mathrm{Spec}(G) \setminus \{0\}$.

Suppose that $\theta^{\pm }$ is an integer and $\theta^{\pm } \in \Phi_{(a_1,0)(a_2,0)}$ for some $\theta \in \mathrm{Spec}(G) \setminus \{0\}$. Then $\theta$ is a positive integer and $(rn_2+\theta-2)^2+4n_2(2r-\theta)$ is a perfect square. Also, 
$4n_2(2r-\theta)$ is an even positive integer. So $(rn_2+\theta-2)^2$ and $(rn_2+\theta-2)^2+4n_2(2r-\theta)$ must have the same parity. Therefore, we have
$$(rn_2+\theta-2)^2+4n_2(2r-\theta) = (rn_2+\theta+2\alpha)^2,~\mathrm{where}~\alpha \in \mathbb{N} \cup \{0\}.$$
This implies that 
$$\theta = \frac{(1-\alpha)(1+rn_2+\alpha)}{(1+n_2+\alpha)}.$$
If $\alpha \in \mathbb{N}$, then $\theta \leq 0$, a contradiction. If $\alpha=0$, then  
$\theta = r - (r - 1)/(1+n_2)$. This implies $(r - 1)/(1+n_2)$ is an integer, which is again a contradiction. Thus $\theta^{\pm } \notin \Phi_{(a_1,0)(a_2,0)}$ for each $\theta \in \mathrm{Spec}(G) \setminus \{0\}$. This implies that $\Phi_{(a_1,0)(a_2,0)}$ is an empty set, a contradiction.

Suppose that $\theta^{\pm }$ is a quadratic integer and $\theta^{\pm } \in \Phi_{(a_1,0)(a_2,0)}$ for some $\theta \in \mathrm{Spec}(G) \setminus \{0\}$. Then Theorem~\ref{Conditions for the existence of Pair-LPST} yields that $\theta^{\pm }=\frac{1}{2}(u+v_{\pm}\sqrt{\Delta})$, where $\Delta$ is a square-free integer such that $\Delta \geq 2$, $u,v_{\pm} \in \mathbb{Z}$ and, $v_{+}$ and $v_{-}$ have the same parity. This, along with Theorem~\ref{Non-bipartite Spectrum}, give 
\begin{equation}\label{b_1}
u + \frac{1}{2}(v_{+} + v_{-})\sqrt{\Delta} = \theta^{+} + \theta^{-} = rn_2+\theta+2.
\end{equation}
We have
\begin{equation}\label{x_9}
4\theta^{+}\theta^{-} = u^2 + v_{+}v_{-}\Delta + u(v_{+}+v_{-})\sqrt{\Delta}.
\end{equation}
Also, From Theorem~\ref{Non-bipartite Spectrum} we have
\begin{equation}\label{d_1}
4\theta^{+}\theta^{-} = 4\theta(n_2+2).
\end{equation}
Now Equation~(\ref{b_1}), Equation~(\ref{x_9}) and Equation~(\ref{d_1}) altogether yields that 
\begin{equation}\label{d_2}
u^2 + v_{+}v_{-}\Delta + u(v_{+}+v_{-})\sqrt{\Delta} = 4(n_2+2)(u-rn_2-2)+2(n_2+2)(v_{+}+v_{-})\sqrt{\Delta}. 
\end{equation}
As $\sqrt{\Delta}$ is an irrational number, Equation~(\ref{d_2}) yields that 
\begin{equation}\label{x_{10}}
(v_{+}+v_{-})(u-2(n_2+2))=0.
\end{equation}
If $v_{+}+v_{-} = 0$, then from Equation~(\ref{b_1}) we obtain that $\theta = u - rn_2 - 2$. This implies that 
$$\Phi_{(a_1,0)(a_2,0)} = \left\{{\theta}^{+},{\theta}^{-}\right\}.$$
Therefore by Lemma~\ref{Laplacian perfect pair state transfer lemma non-bipartite}, we have 
$\Phi^{+}_{(a_1,0)(a_2,0),(b_1,0)(b_2,0)} = \{{\theta}^{+},{\theta}^{-}\}$ or $\Phi^{-}_{(a_1,0)(a_2,0),(b_1,0)(b_2,0)} = \{{\theta}^{+},{\theta}^{-}\}$.
This implies that either $\Phi^{-}_{(a_1,0)(a_2,0),(b_1,0)(b_2,0)}$ is empty or $\Phi^{+}_{(a_1,0)(a_2,0),(b_1,0)(b_2,0)}$ is empty, a contradiction. If $v_{+}+v_{-} \neq 0$, then from Equation~(\ref{b_1}) we obtain that $u = 2n_2+4$. Now we have 
\begin{align*}
&\frac{1}{2}(v_{+}+v_{-})\sqrt{\Delta} = (r-2)n_2+\theta-2~~\mathrm{and}\\
&\frac{1}{2}(v_{+}-v_{-})\sqrt{\Delta} = \sqrt{(rn_2+\theta-2)^2+4n_2(2r-\theta)}.
\end{align*}
The preceding two equations imply that
\begin{equation}\label{d_5}
v_{+}v_{-}\Delta = 4n_2(n_2+2)(1-r).
\end{equation}
Recall that $v_{+}$ and $v_{-}$ are of the same parity. If both of $v_{+}$ and $v_{-}$ are odd. Then Equation~(\ref{d_5}) implies that $4 | \Delta$, contradicts the fact that $\Delta$ is a square-free integer. If both of $v_{+}$ and $v_{-}$ are even, then we write $v_{+} = 2\alpha$ and $v_{-} = 2\beta$. This implies that $\theta = (2-r)n_2 + 2 + (\alpha+\beta)\sqrt{\Delta}$. Let $\theta^\prime = (2-r)n_2 + 2 - (\alpha+\beta)\sqrt{\Delta}$. Since the characteristic polynomial of $L_{G}$ has integer coefficients, the algebraic conjugate of $\theta$, that is, $\theta^\prime$ is also an eigenvalue of $L_{G}$. Since both of $\theta$ and $\theta^\prime$ are positive, we have $\theta + \theta^\prime > 0$. This implies that $(r-2)n_2 = 1$. Since $r\geq3$, we have that $r=3$ and $n_2=1$. Now we have
$$\theta = 1 + (\alpha+\beta)\sqrt{\Delta}~~\mathrm{and}~~\theta^\prime = 1 - (\alpha+\beta)\sqrt{\Delta}.$$
This implies that $-1 < (\alpha+\beta)\sqrt{\Delta} < 1$, that is, $|\alpha+\beta|\sqrt{\Delta} < 1$. This is a contradiction as $\Delta \geq 2$ and $\alpha+\beta$ is a non-zero integer.

Now we consider that $G$ is bipartite and $G \diamond \mathcal{H}$ admits Pair-LPST from $((a_1,0),(a_2,0))$ to $((b_1,0),(b_2,0))$. We have $\Phi_{(a_1,0)(a_2,0)} \subseteq \{\theta^{\pm } \colon \theta \in \mathrm{Spec}(G)\setminus \{0,2r\}\} \cup \{r(n_2+2)\}$. Let $V_1 \cup V_2$ be the bipartition of $G$. Then the following two cases arise.

\noindent \textbf{Case 1.} Let $a_1, a_2 \in V_1$ or $a_1, a_2 \in V_2$. Then Theorem~\ref{Bipartite Spectrum} yields that $r(n_2+2) \notin \Phi_{(a_1,0)(a_2,0)}$. This implies $\Phi_{(a_1,0)(a_2,0)} \subseteq \{\theta^{\pm } \colon \theta \in \mathrm{Spec}(G)\setminus \{0,2r\}\}$. Now following the procedure of the non-bipartite case, we arrive at a contradiction. 

\noindent \textbf{Case 2.} Let $a_1 \in V_1,~a_2 \in V_2$ or $a_2 \in V_1,~a_1 \in V_2$. Then Theorem~\ref{Bipartite Spectrum} yields that $r(n_2+2) \in \Phi_{(a_1,0)(a_2,0)}$. Since $((a_1,0),(a_2,0))$ involves in Pair-LPST, Theorem~\ref{Conditions for the existence of Pair-LPST} implies that all the eigenvalues in $\Phi_{(a_1,0)(a_2,0)}$ are integers. Now one can prove that $\theta^{\pm } \notin \Phi_{(a_1,0)(a_2,0)}$ for all $\theta \in \mathrm{Spec}(G)\setminus \{0,2r\}$. This implies that $\Phi_{(a_1,0)(a_2,0)} = \{r(n_2+2)\}$, which is a contradiction to Lemma~\ref{fixed state}. Therefore, the pair $((a_1,0),(a_2,0))$ cannot be involved in Pair-LPST. This completes the proof.
\end{proof}

\begin{theorem}\label{non-existence of Pair-LPST 2nd Lemma}
Let $a_1$ be a vertex of $G$, $e_s$ be an edge of $G$ and $a_2$ be a vertex of $H_s$. Further, let $(r-1)/(1+n_2)$ not be an integer. Then the pair $((a_1,0),(a_2,e_s))$ cannot be involved in Pair-LPST.  
\end{theorem}
\begin{proof}
First consider that $G$ is non-bipartite. If possible, let $((a_1,0),(a_2,e_s))$ be involved in Pair-LPST. Then From the proof of Lemma~\ref{Laplacian perfect pair state transfer lemma non-bipartite}, we obtain that $G \diamond \mathcal{H}$ admits Pair-LPST from $((a_1,0),(a_2,e_s))$ to $((b_1,0),(b_2,e_t))$ for some vertex $b_1$ of $G$, edge $e_t$ of $G$ and vertex $b_2$ of $H_t$. From Equation~(\ref{x_3}), it is clear that $rn_2+2 \in \Phi_{(a_1,0)(a_2,e_s)}$. Since $((a_1,0),(a_2,e_s))$ involves in Pair-LPST, Theorem~\ref{Conditions for the existence of Pair-LPST} implies that all the eigenvalues in $\Phi_{(a_1,0)(a_2,e_s)}$ are integers. From Theorem~\ref{Non-bipartite Spectrum}, we recall that 
\begin{equation}\label{a_1}
\theta^{\pm }=\frac{rn_2+\theta+2 \pm \sqrt{(rn_2+\theta-2)^2+4n_2(2r-\theta)}}{2}.
\end{equation} 
If possible, let $\theta^{+} \in \Phi_{(a_1,0)(a_2,e_s)}$ for some $\theta \in \mathrm{Spec}(G) \setminus \{0\}$. Now we have 
$$rn_2+\theta+2 + \sqrt{(rn_2+\theta-2)^2+4n_2(2r-\theta)}=2 \theta^{+}.$$
This implies that 
$${(rn_2+\theta-2)}^2+4n_2(2r-\theta)={(2 \theta^{+}-4)}^2+{(rn_2+\theta-2)}^2-2(2 \theta^{+}-4)(rn_2+\theta-2).$$
This further implies that 
\begin{equation}\label{a_2}
(\theta^{+}-2-n_2)\theta={(\theta^{+}-2)}^2-(\theta^{+}-2)(rn_2-2)-2rn_2.
\end{equation} 
If $\theta^{+}-2-n_2=0$, then we obtain that $r=1$, contradicts the fact that $r\geq2$. Therefore, $\theta^{+}-2-n_2$ is a non-zero integer. Then Equation~(\ref{a_2}) implies that $\theta$ is an integer. But $\theta^{+} + \theta^{-} = rn_2+\theta+2$. Therefore, we obtain that $\theta^{-}$ is also an integer. In a similar way, it can be proved that if $\theta^{-} \in \Phi_{(a_1,0)(a_2,e_s)}$ for some $\theta \in \mathrm{Spec}(G) \setminus \{0\}$, then $\theta^{+}$ is an integer. Now proceeding as in the proof of Theorem~\ref{non-existence of Pair-LPST First Theorem}, we arrive at a contradiction. Therefore, we have 
$$\Phi_{(a_1,0)(a_2,e_s)} \subseteq \{rn_2+2, \mu+2 \colon \mu \in \mathrm{Spec}(H_j) ~\mathrm{for} ~1 \leq j \leq m\}.$$
Consider that 
$$\Phi_{(a_1,0)(a_2,e_s)} = \{rn_2+2, \mu_{j_1}+2,\ldots,\mu_{j_k}+2\}.$$
Now post-multiplying by $\textbf{e}_{(a_1,0)} - \textbf{e}_{(a_2,e_s)}$ on the both sides of
$$\sum_{\theta \in \mathrm{Spec}(G) \setminus \{0\}}F_{\theta^{\pm}}(G \diamond \mathcal{H}) + F_{0^{\pm}}(G \diamond \mathcal{H}) + \sum_{\mu \in \cup_{j=1}^{m}\mathrm{Spec}(H_j)}F_{\mu}(G \diamond \mathcal{H}) = I_{(n_1+mn_2)\times(n_1+mn_2)},$$
we have
$$\left(F_{rn_2+2}(G \diamond \mathcal{H}) + \sum_{i=1}^{k}F_{\mu_{j_i}+2}(G \diamond \mathcal{H})\right)\left(\textbf{e}_{(a_1,0)} - \textbf{e}_{(a_2,e_s)}\right) = \textbf{e}_{(a_1,0)} - \textbf{e}_{(a_2,e_s)}.$$
Now pre-multiplying both sides of the preceding equation by $\textbf{e}^T_{(a_1,0)}$ and then simplifying, we obtain that
$$\frac{rn_2}{n_1(rn_2+2)}+\frac{2}{n_1(rn_2+2)}=1.$$
This implies that $n_1=1$, a contradiction. Therefore, the pair $((a_1,0),(a_2,e_s))$ cannot be involved in Pair-LPST.  

Now we consider the case that $G$ is bipartite. Proceeding as in the non-bipartite case, we obtain that $\theta^{\pm} \notin \Phi_{(a_1,0)(a_2,e_s)}$ for all $\theta \in \mathrm{Spec}(G) \setminus \{0\}$. Therefore, 
$$\Phi_{(a_1,0)(a_2,e_s)} \subseteq \{rn_2+2,r(n_2+2), \mu+2 \colon \mu \in \mathrm{Spec}(H_j) ~\mathrm{for} ~1 \leq j \leq m\}.$$
The following two cases arise.

\noindent \textbf{Case 1.} If $r(n_2+2) \notin \Phi_{(a_1,0)(a_2,e_s)}$, then 
$$\Phi_{(a_1,0)(a_2,e_s)} \subseteq \{rn_2+2,\mu+2 \colon \mu \in \mathrm{Spec}(H_j) ~\mathrm{for} ~1 \leq j \leq m\}.$$
The rest of the proof is similar to that of the non-bipartite case.

\noindent \textbf{Case 2.} If $r(n_2+2) \in \Phi_{(a_1,0)(a_2,e_s)}$, then the following two sub-cases arise. 
 
\noindent \textbf{Sub-case 1.} $\Phi_{(a_1,0)(a_2,e_s)} = \{rn_2+2,r(n_2+2)\}$. Now post-multiplying by $\textbf{e}_{(a_1,0)} - \textbf{e}_{(a_2,e_s)}$ on the both sides of
\begin{align*}
&\sum_{\theta \in \mathrm{Spec}(G) \setminus \{0,2r\}}F_{\theta^{\pm}}(G \diamond \mathcal{H}) + F_{0^{\pm}}(G \diamond \mathcal{H}) + F_{r(n_2+2)}(G \diamond \mathcal{H}) + \sum_{\mu \in \cup_{j=1}^{m}\mathrm{Spec}(H_j)}F_{\mu}(G \diamond \mathcal{H}) \\
&= I_{(n_1+mn_2)\times(n_1+mn_2)},
\end{align*}
we have
$$\left(F_{rn_2+2}(G \diamond \mathcal{H}) + F_{r(n_2+2)}(G \diamond \mathcal{H})\right)\left(\textbf{e}_{(a_1,0)} - \textbf{e}_{(a_2,e_s)}\right) = \textbf{e}_{(a_1,0)} - \textbf{e}_{(a_2,e_s)}.$$ 
Now pre-multiplying both sides of the preceding equation by $\textbf{e}^T_{(a_1,0)}$ and then simplifying, we obtain that
$$\frac{rn_2+2}{n_1(rn_2+2)}+\frac{1}{n_1}=1.$$
This implies that $n_1=2$, a contradiction. 

\noindent \textbf{Sub-case 2.} $\Phi_{(a_1,0)(a_2,e_s)} = \{rn_2+2,r(n_2+2), \mu_{j_1}+2,\ldots,\mu_{j_k}+2\}$. Note that for all $\mu_{j_i}$, we have
$$\textbf{e}^T_{(a_1,0)}F_{\mu_{j_i}+2}(G \diamond \mathcal{H})\textbf{e}_{(a_1,0)} = \textbf{0}_{n_1+mn_2} = \textbf{e}^T_{(a_1,0)}F_{\mu_{j_i}+2}(G \diamond \mathcal{H})\textbf{e}_{(a_2,e_s)}.$$
Now the rest of the proof is similar to that of the Sub-case 1. Therefore, the pair $((a_1,0),(a_2,e_s))$ cannot be involved in Pair-LPST. This completes the proof.
\end{proof}

\section{Pair-LPGST on edge coronas}\label{Pair-LPGST on edge corona}
In this section, we explore the existence of Pair-LPGST on the edge corona $G \diamond \mathcal{H}$. Unless otherwise stated, $G$ is a connected $r$-regular graph with $r \geq 2$, $|V(G)|=n_1$ and $|E(G)|=m$. Let $\mathcal{H}=(H_1,\ldots,H_m)$, where $H_j$ is connected and $|V(H_j)|=n_2$ for $1 \leq j \leq m$. Also, let $\Delta_\theta = \sqrt{(r n_2 + \theta - 2)^2+4n_2(2r-\theta)}$, where $\theta \in \mathrm{Spec}(G)$ if $G$ is non-bipartite and $\theta \in \mathrm{Spec}(G) \setminus \{2r\}$ if $G$ is bipartite. If $H_1 = \cdots = H_m = H$, then the next theorem follows from Monterde and Pal~\cite{MonterdePal2026} [see Theorem 7].
\begin{theorem}\emph{\cite{MonterdePal2026}}\label{Graphs with clusters Pair-LPGST} 
Let $e_s$ be a fixed edge of $G$ and $h_i$ be a vertex of $H$, where $1 \leq i \leq 4$. Then the graph $H$ admits Pair-LPGST from $(h_1,h_2)$ to $(h_3,h_4)$
if and only if the edge corona $G \diamond H$ admits Pair-LPGST from $((h_1,e_s),(h_2,e_s))$ to $((h_3,e_s),(h_4,e_s))$.
\end{theorem}
Theorem~\ref{Graphs with clusters Pair-LPGST} fails to give any information about the existence of Pair-LPGST from pair $((a_1,0), (a_2,0))$ to $((b_1,0), (b_2,0))$, where $a_1,a_2,b_1,b_2 \in V(G)$. This motivates us to explore  the existence of Pair-LPGST on $G \diamond \mathcal{H}$ from pair of the form $((a_1,0), (a_2,0))$ to pair of the form $((b_1,0), (b_2,0))$. We begin with the following theorem of Li et al.~\cite{Edgecoronastatetransfer2022}. 

\begin{theorem}\emph{\cite{Edgecoronastatetransfer2022}}\label{Edge corona entry}
Let $a_1$ and $a_2$ be two vertices of the graph $G$. If $G$ is non-bipartite, then
\begin{align*}
&\textbf{e}^T_{(a_1,0)}U_{G \diamond \mathcal{H}}(t)\textbf{e}_{(a_2,0)}\\
&=\sum_{\theta \in \mathrm{Spec}(G)}\exp(-\textbf{i}t(rn_2 + \theta + 2/2))F_{\theta}(G)_{a_1a_2}\left(\cos\left(\frac{\Delta_\theta t}{2}\right) + \textbf{i}\frac{2-r n_2 - \theta}{\Delta_\theta}\sin\left(\frac{\Delta_\theta t}{2}\right)\right).
\end{align*} 
If $G$ is bipartite, then
\begin{align*}
&\textbf{e}^T_{(a_1,0)}U_{G \diamond \mathcal{H}}(t)\textbf{e}_{(a_2,0)}\\
&=\sum_{\theta \in \mathrm{Spec}(G) \setminus\{2r\}}\exp(-\textbf{i}t(rn_2 + \theta + 2/2))F_{\theta}(G)_{a_1a_2}\left(\cos\left(\frac{\Delta_\theta t}{2}\right) + \textbf{i}\frac{2-r n_2 - \theta}{\Delta_\theta}\sin\left(\frac{\Delta_\theta t}{2}\right)\right)\\\
&+ \exp(-\textbf{i}t(r n_2+2))F_{2r}(G)_{a_1a_2}.
\end{align*}
\end{theorem}

The next lemma is an immediate consequence of Theorem~\ref{Edge corona entry}; however, we supply the proof for completeness. 
\begin{lema}\label{Edge corona lemma}
Let $a_1, a_2, b_1$ and $b_2$ be four vertices of the graph $G$. If $G$ is non-bipartite, then
\begin{align*}
&\frac{1}{2}{\left(\textbf{e}_{(a_1,0)} - \textbf{e}_{(a_2,0)}\right)}^TU_{G \diamond \mathcal{H}}(t)\left(\textbf{e}_{(b_1,0)} - \textbf{e}_{(b_2,0)}\right)\\
&= \frac{1}{2}\sum_{\theta \in \mathrm{Spec}(G)}\exp(-\textbf{i}t((r n_2 + \theta + 2)/2)){\left(\textbf{e}_{a_1} - \textbf{e}_{a_2}\right)}^TF_{\theta}(G)\left(\textbf{e}_{b_1} - \textbf{e}_{b_2}\right)\\
&\times \left(\cos\left(\frac{\Delta_\theta t}{2}\right) + \textbf{i}\frac{2-r n_2 - \theta}{\Delta_\theta}\sin\left(\frac{\Delta_\theta t}{2}\right)\right).
\end{align*}
If $G$ is bipartite, then
\begin{align*}
&\frac{1}{2}{\left(\textbf{e}_{(a_1,0)} - \textbf{e}_{(a_2,0)}\right)}^TU_{G \diamond \mathcal{H}}(t)\left(\textbf{e}_{(b_1,0)} - \textbf{e}_{(b_2,0)}\right)\\
&= \frac{1}{2}\sum_{\theta \in \mathrm{Spec}(G) \setminus \{2r\}}\exp(-\textbf{i}t((r n_2 + \theta + 2)/2)){\left(\textbf{e}_{a_1} - \textbf{e}_{a_2}\right)}^TF_{\theta}(G)\left(\textbf{e}_{b_1} - \textbf{e}_{b_2}\right)\\
& \times \left(\cos\left(\frac{\Delta_\theta t}{2}\right) + \textbf{i}\frac{2-r n_2 - \theta}{\Delta_\theta}\sin\left(\frac{\Delta_\theta t}{2}\right)\right)
+ \frac{1}{2}\exp(-\textbf{i}t(r n_2+2)) {(\textbf{e}_{a_1} - \textbf{e}_{a_2})}^TF_{2r}(G)(\textbf{e}_{b_1} - \textbf{e}_{b_2}).
\end{align*}
\end{lema}
\begin{proof}
If $G$ is non-bipartite, then applying Theorem~\ref{Edge corona entry} we obtain that
\begin{align*}
&\frac{1}{2}{\left(\textbf{e}_{(a_1,0)} - \textbf{e}_{(a_2,0)}\right)}^TU_{G \diamond \mathcal{H}}(t)\left(\textbf{e}_{(b_1,0)} - \textbf{e}_{(b_2,0)}\right)\\
&=\frac{1}{2}\left(\textbf{e}^T_{(a_1,0)}U_{G \diamond \mathcal{H}}(t)\textbf{e}_{(b_1,0)} - \textbf{e}^T_{(a_1,0)}U_{G \diamond \mathcal{H}}(t)\textbf{e}_{(b_2,0)} - \textbf{e}^T_{(a_2,0)}U_{G \diamond \mathcal{H}}(t)\textbf{e}_{(b_1,0)} + \textbf{e}^T_{(a_2,0)}U_{G \diamond \mathcal{H}}(t)\textbf{e}_{(b_2,0)}\right)\\
&= \frac{1}{2}\sum_{\theta \in \mathrm{Spec}(G)}\exp(-\textbf{i}t((r n_2 + \theta + 2)/2)){\left(\textbf{e}_{a_1} - \textbf{e}_{a_2}\right)}^TF_{\theta}(G)\left(\textbf{e}_{b_1} - \textbf{e}_{b_2}\right)\\
&\times \left(\cos\left(\frac{\Delta_\theta t}{2}\right) + \textbf{i}\frac{2-r n_2 - \theta}{\Delta_\theta}\sin\left(\frac{\Delta_\theta t}{2}\right)\right).
\end{align*}
If $G$ is bipartite, then the proof is similar to that of the preceding part, and hence the details are omitted to avoid repetitive argument. This completes the proof.
\end{proof}

In the following, we present a lemma which is important for the investigation of Pair-LPGST on $G \diamond \mathcal{H}$. 
\begin{lema}\label{Preliminary Lemma}
Let $a_1$ and $a_2$ be two vertices of $G$. Also, let $\theta \in \Phi_{a_1a_2}$ such that $\theta$ is an integer and $\Delta_{\theta}$ is rational. Then $(r-1)/(1+n_2)$ must be an integer.
\end{lema}
\begin{proof}
Recall that $\theta^{\pm}=\frac{rn_2+\theta+2\pm \Delta_{\theta}}{2}$. Since $\theta \in \Phi_{a_1a_2}$ such that $\theta$ is an integer and $\Delta_{\theta}$ is rational, we obtain that both of $\theta^{+}$ and $\theta^{-}$ are integers. Now following the proof of Theorem~\ref{non-existence of Pair-LPST First Theorem}, we obtain that $\theta = r - (r-1)/(1+n_2)$. This implies that $(r-1)/(1+n_2)$ is an integer. 
\end{proof}

Now we present two theorems that give sufficient conditions for the existence of Pair-LPGST on edge coronas. The following theorem tells that if $G$ is non-bipartite, then under some mild conditions the edge corona $G \diamond \mathcal{H}$ admits Pair-LPGST. 

\begin{theorem}\label{Edge Corona First Theorem PGST}
Let $G$ be non-bipartite and admit Pair-LPST from $(a_1, a_2)$ to $(b_1, b_2)$ at time $\frac{\pi}{2}$. Also, let $n_2$ be an even integer and $(r-1)/(1+n_2)$ not be an integer. Then the edge corona $G \diamond \mathcal{H}$ admits Pair-LPGST from $((a_1,0), (a_2,0))$ to $((b_1,0), (b_2,0))$.  
\end{theorem}
\begin{proof}
It is clear that $0 \notin \Phi_{a_1a_2}$. Since $G$ admits Pair-LPST from $(a_1, a_2)$ to $(b_1, b_2)$ at time $\frac{\pi}{2}$, applying Theorem~\ref{Conditions for the existence of Pair-LPST} and Corollary~\ref{minimum time Corollary}, we obtain that all the eigenvalues in $\Phi_{a_1a_2}$ are positive integers. We have $\Delta_\theta = \sqrt{(r n_2 + \theta - 2)^2+4n_2(2r-\theta)}$, where $\theta \in \mathrm{Spec}(G)$. This implies that $\Delta^2_{\theta}$ is an integer for all $\theta \in \Phi_{a_1a_2}$. If $\theta \in \Phi_{a_1a_2}$, then we write $\Delta_{\theta} = u_{\theta} \sqrt{v_{\theta}}$, where $u_{\theta}$ is an integer and $v_{\theta}$ is a square free integer. Since $(r-1)/(1+n_2)$ is not an integer, it follows from Lemma~\ref{Preliminary Lemma} that $\Delta_{\theta}$ is irrational for all $\theta \in \Phi_{a_1a_2}$. This implies that $v_{\theta}$ is strictly greater than 1 for all $\theta \in \Phi_{a_1a_2}$. Therefore, Lemma~\ref{Linearly independent corollary} implies that the set $\{1\}\cup \{\sqrt{v_{\theta}} \colon \theta \in \Phi_{a_1a_2}\}$ is linearly independent over the field of rational numbers. 

Let $\varepsilon$ be an arbitrary positive real number. Since the cosine function is uniformly continuous on $\mathbb{R}$, there exists a positive real number $\delta$ such that $|x_1-x_2|< \delta$ implies $|\cos x_1 - \cos x_2|< \frac{\varepsilon}{2}$ for all $x_1, x_2 \in \mathbb{R}$. Let max$\{u_{\theta} \colon \theta \in \Phi_{a_1a_2}\}=u_{\theta_i}$. Since the set $\{1\}\cup \{\sqrt{v_{\theta}} \colon \theta \in \Phi_{a_1a_2}\}$ is linearly independent over the field of rational numbers, applying Lemma~\ref{Kronecker approximation theorem} yields that there exists $\ell, s_{\theta} \in \mathbb{Z}$ such that $\left|\ell \sqrt{v_\theta} + \frac{\sqrt{v_\theta}}{4} - s_\theta\right|<\frac{2 \delta}{4 u_{\theta_i}\pi}$ for all $\theta \in \Phi_{a_1a_2}$. This implies that $\left|\ell \sqrt{v_\theta} + \frac{\sqrt{v_\theta}}{4} - s_\theta\right|<\frac{2 \delta}{4 u_{\theta}\pi}$ for all $\theta \in \Phi_{a_1a_2}$. This further implies that $\left|(4\ell+1)\frac{\pi \Delta_\theta}{2} - 2u_{\theta}s_{\theta}\pi \right|<\delta$ for all $\theta \in \Phi_{a_1a_2}$. Putting $t=(4\ell+1)\pi$, we obtain that $\left|\frac{\Delta_\theta t}{2} - 2u_{\theta}s_{\theta}\pi \right|<\delta$ for all $\theta \in \Phi_{a_1a_2}$. This implies that $\left|\cos \left(\frac{\Delta_\theta t}{2}\right) - \cos(2u_\theta s_\theta \pi)\right|<\frac{\varepsilon}{2}$ for all $\theta \in \Phi_{a_1a_2}$. Therefore, $\left|\cos \left(\frac{\Delta_\theta t}{2}\right) - 1 \right|<\frac{\varepsilon}{2}$ for all $\theta \in \Phi_{a_1a_2}$. Similarly, we obtain that 
$\left|\sin \left(\frac{\Delta_\theta t}{2}\right)\right|<\frac{\varepsilon}{2}$ for all $\theta \in \Phi_{a_1a_2}$. It is clear that $\left|\frac{2-rn_2-\theta}{\Delta_\theta}    \right|<1$ for all $\theta \in \Phi_{a_1a_2}$. Therefore, $\left|\cos\left(\frac{\Delta_\theta t}{2}\right) + \textbf{i}\frac{2-rn_2 - \theta}{\Delta_\theta}\sin\left(\frac{\Delta_\theta t}{2}\right)-1 \right|<\varepsilon$ for all $\theta \in \Phi_{a_1a_2}$. Since $n_2$ is an even integer, $\theta \in \Phi_{a_1a_2}$ are integers and $t=(4\ell+1)\pi$, applying Lemma~\ref{Edge corona lemma} we obtain that 
\begin{align*}
&\frac{1}{2}{\left(\textbf{e}_{(a_1,0)} - \textbf{e}_{(a_2,0)}\right)}^TU_{G \diamond \mathcal{H}}(t)\left(\textbf{e}_{(b_1,0)} - \textbf{e}_{(b_2,0)}\right)\\
&= \sum_{\theta \in \Phi_{a_1a_2}}\frac{(-1)^{\frac{rn_2+2}{2}}}{2}\exp\left(\frac{-\textbf{i}\pi \theta}{2}\right)
{\left(\textbf{e}_{a_1} - \textbf{e}_{a_2}\right)}^TF_{\theta}(G)\left(\textbf{e}_{b_1} - \textbf{e}_{b_2}\right)\\
&\times \left(\cos\left(\frac{\Delta_\theta t}{2}\right) + \textbf{i}\frac{2-rn_2 - \theta}{\Delta_\theta}\sin\left(\frac{\Delta_\theta t}{2}\right)\right).
\end{align*}
It can be verified that 
$$\sum_{\theta \in \Phi_{a_1a_2}}\left|\frac{1}{2}\exp\left(\frac{-\textbf{i}\pi \theta}{2}\right){({\textbf{e}_{a_1}} - {\textbf{e}_{a_2}})}^T
F_{\theta}(G)({\textbf{e}_{b_1}} - {\textbf{e}_{b_2}})\right|=1.$$
Now we have 
\begin{align*}
&\frac{1}{2}{\left(\textbf{e}_{(a_1,0)} - \textbf{e}_{(a_2,0)}\right)}^TU_{G \diamond \mathcal{H}}(t)\left(\textbf{e}_{(b_1,0)} - \textbf{e}_{(b_2,0)}\right) - \frac{(-1)^{\frac{rn_2+2}{2}}}{2}{\left(\textbf{e}_{a_1} - \textbf{e}_{a_2}\right)}^TU_G\left(\frac{\pi}{2}\right)\left(\textbf{e}_{b_1} - \textbf{e}_{b_2}\right)\\
&=\sum_{\theta \in \Phi_{a_1a_2}}\frac{(-1)^{\frac{rn_2+2}{2}}}{2}\exp\left(\frac{-\textbf{i}\pi \theta}{2}\right){\left(\textbf{e}_{a_1} - \textbf{e}_{a_2}\right)}^T F_{\theta}(G)\left(\textbf{e}_{b_1}-\textbf{e}_{b_2}\right)\\
&\times \left(\cos\left(\frac{\Delta_\theta t}{2}\right) + \textbf{i}\frac{2-rn_2 - \theta}{\Delta_\theta}\sin\left(\frac{\Delta_\theta t}{2}\right)-1\right).
\end{align*}
This implies 
$$\left| \frac{1}{2}{\left(\textbf{e}_{(a_1,0)} - \textbf{e}_{(a_2,0)}\right)}^TU_{G \diamond \mathcal{H}}(t)\left(\textbf{e}_{(b_1,0)} - \textbf{e}_{(b_2,0)}\right) - \frac{(-1)^{\frac{rn_2+2}{2}}}{2}{\left(\textbf{e}_{a_1} - \textbf{e}_{a_2}\right)}^TU_G\left(\frac{\pi}{2}\right)\left(\textbf{e}_{b_1} - \textbf{e}_{b_2}\right) \right|<\varepsilon.$$
Since $G$ admits Pair-LPST from $(a_1, a_2)$ to $(b_1, b_2)$ at time $\frac{\pi}{2}$, we have 
$$\left|\frac{1}{2}{({\textbf{e}_{a_1}} - {\textbf{e}_{a_2}})}^TU_G\left(\frac{\pi}{2}\right)({\textbf{e}_{b_1}} - {\textbf{e}_{b_2}})\right| = 1.$$
Now we obtain that
$$\left|\left|\frac{1}{2}{\left(\textbf{e}_{(a_1,0)} - \textbf{e}_{(a_2,0)}\right)}^TU_{G \diamond \mathcal{H}}(t)\left(\textbf{e}_{(b_1,0)} - \textbf{e}_{(b_2,0)}\right)\right|-1\right|<\varepsilon.$$
Therefore, the edge corona $G \diamond \mathcal{H}$ admits Pair-LPGST from $((a_1,0), (a_2,0))$ to $((b_1,0), (b_2,0))$. 
\end{proof}

Let $G$ be an $r$-regular bipartite with the bipartition $V_1 \cup V_2$. Then $2r$ is an eigenvalue of $G$ and the corresponding eigenprojector is given by
\begin{equation}\label{H_1}
F_{2r}(G) = \frac{1}{n}\begin{pmatrix}
         J_{|V_1|\times |V_1|}  &  -J_{|V_1|\times |V_2|}\\
         -J_{|V_2|\times |V_1|} &  J_{|V_2|\times |V_2|}
         \end{pmatrix}.
\end{equation}

The next theorem tells that if $G$ is bipartite, then under some mild conditions the edge corona $G \diamond \mathcal{H}$ admits Pair-LPGST.

\begin{theorem}\label{Edge Corona Second Theorem PGST}
Let $G$ be bipartite with the bipartition $V_1 \cup V_2$. Also, let $G$ admit Pair-LPST from $(a_1, a_2)$ to $(b_1, b_2)$ at time $\frac{\pi}{2}$, $n_2$ be an even integer and $(r-1)/(1+n_2)$ not be an integer.
\begin{enumerate}
\item[(a)] If $a_1, a_2 \in V_1$ or $a_1, a_2 \in V_2$, then the edge corona $G \diamond \mathcal{H}$ admits Pair-LPGST from $((a_1,0), (a_2,0))$ to $((b_1,0), (b_2,0))$. 
\item[(b)] If $a_1 \in V_1,~a_2 \in V_2$ or $a_2 \in V_1,~a_1 \in V_2$ and $r(n_2-2) \equiv 2~(\mathrm{mod}~4)$, then the edge corona $G \diamond \mathcal{H}$ admits Pair-LPGST from $((a_1,0), (a_2,0))$ to $((b_1,0), (b_2,0))$. 
\end{enumerate}
\end{theorem}
\begin{proof}
Let $\varepsilon$ be an arbitrary positive real number. Following the notations and the proof of Theorem~\ref{Edge Corona First Theorem PGST}, we obtain that
$$\left|\cos\left(\frac{\Delta_\theta t}{2}\right) + \textbf{i}\frac{2-rn_2 - \theta}{\Delta_\theta}\sin\left(\frac{\Delta_\theta t}{2}\right)-1 \right|<\varepsilon$$ 
for all $\theta \in \Phi_{a_1a_2} \setminus \{2r\}$, where $t=(4\ell+1)\pi$.

\noindent $(a)$  If $a_1, a_2 \in V_1$ or $a_1, a_2 \in V_2$, then from Equation~(\ref{H_1}) it is clear that $2r$ does not belong to $\Phi_{a_1a_2}$. From Lemma~\ref{Edge corona lemma} and proceeding as in the proof of Theorem~\ref{Edge Corona First Theorem PGST}, we have 
$$\left| \frac{1}{2}{\left(\textbf{e}_{(a_1,0)} - \textbf{e}_{(a_2,0)}\right)}^TU_{G \diamond \mathcal{H}}(t)\left(\textbf{e}_{(b_1,0)} - \textbf{e}_{(b_2,0)}\right) - \frac{(-1)^{\frac{rn_2+2}{2}}}{2}{\left(\textbf{e}_{a_1} - \textbf{e}_{a_2}\right)}^TU_G\left(\frac{\pi}{2}\right)\left(\textbf{e}_{b_1}-\textbf{e}_{b_2}\right) \right|<\varepsilon.$$
This implies that
$$\left|\left|\frac{1}{2}{\left(\textbf{e}_{(a_1,0)} - \textbf{e}_{(a_2,0)}\right)}^TU_{G \diamond \mathcal{H}}(t)\left(\textbf{e}_{(b_1,0)} - \textbf{e}_{(b_2,0)}\right)\right|-1\right|<\varepsilon.$$
Therefore, the edge corona $G \diamond \mathcal{H}$ admits Pair-LPGST from $((a_1,0), (a_2,0))$ to $((b_1,0), (b_2,0))$. 

\noindent $(b)$ If $a_1 \in V_1,~a_2 \in V_2$ or $a_2 \in V_1,~a_1 \in V_2$, then from Equation~(\ref{H_1}) it is clear that $2r$ belongs to $\Phi_{a_1a_2}$. From Lemma~\ref{Edge corona lemma}, we have
\begin{align*}
&\frac{1}{2}{\left(\textbf{e}_{(a_1,0)} - \textbf{e}_{(a_2,0)}\right)}^TU_{G \diamond \mathcal{H}}(t)\left(\textbf{e}_{(b_1,0)} - \textbf{e}_{(b_2,0)}\right)\\
&= \sum_{\theta \in  \Phi_{a_1a_2} \setminus \{2r\}}\frac{1}{2}\exp(-\textbf{i}t((r n_2 + \theta + 2)/2)){\left(\textbf{e}_{a_1} - \textbf{e}_{a_2}\right)}^TF_{\theta}(G)\left(\textbf{e}_{b_1} - \textbf{e}_{b_2}\right)\\
& \times \left(\cos\left(\frac{\Delta_\theta t}{2}\right) + \textbf{i}\frac{2-r n_2 - \theta}{\Delta_\theta}\sin\left(\frac{\Delta_\theta t}{2}\right)\right)
+ \frac{1}{2}\exp(-\textbf{i}t(r n_2+2)) {(\textbf{e}_{a_1} - \textbf{e}_{a_2})}^TF_{2r}(G)(\textbf{e}_{b_1} - \textbf{e}_{b_2}).
\end{align*}
Since for each $\theta \in  \Phi_{a_1a_2}$, $\theta$ is a positive integer, $t=(4\ell+1)\pi$, $n_2$ is an even integer, $r(n_2-2) \equiv 2~(\mathrm{mod}~4)$ and $2r$ is an eigenvalue of $G$, we obtain that 
\begin{align*}
&\sum_{\theta \in  \Phi_{a_1a_2} \setminus \{2r\}}\frac{1}{2}\exp(-\textbf{i}t((r n_2 + \theta + 2)/2)){\left(\textbf{e}_{a_1} - \textbf{e}_{a_2}\right)}^TF_{\theta}(G)\left(\textbf{e}_{b_1} - \textbf{e}_{b_2}\right)\\
&+ \frac{1}{2}\exp(-\textbf{i}t(r n_2+2)) {(\textbf{e}_{a_1} - \textbf{e}_{a_2})}^TF_{2r}(G)(\textbf{e}_{b_1} - \textbf{e}_{b_2})\\
&=\frac{(-1)^{\frac{rn_2+2}{2}}}{2}{\left(\textbf{e}_{a_1} - \textbf{e}_{a_2}\right)}^TU_G\left(\frac{\pi}{2}\right)\left(\textbf{e}_{b_1} - \textbf{e}_{b_2}\right).
\end{align*}
Now we have
\begin{align*}
&\frac{1}{2}{\left(\textbf{e}_{(a_1,0)} - \textbf{e}_{(a_2,0)}\right)}^TU_{G \diamond \mathcal{H}}(t)\left(\textbf{e}_{(b_1,0)} - \textbf{e}_{(b_2,0)}\right) - \frac{(-1)^{\frac{rn_2+2}{2}}}{2}{\left(\textbf{e}_{a_1} - \textbf{e}_{a_2}\right)}^TU_G\left(\frac{\pi}{2}\right)\left(\textbf{e}_{b_1} - \textbf{e}_{b_2}\right)\\
&=\sum_{\theta \in  \Phi_{a_1a_2} \setminus \{2r\}}\frac{{(-1)}^{\frac{rn_2+2}{2}}}{2}
\exp \left(\frac{-\textbf{i} \pi \theta}{2}\right){\left(\textbf{e}_{a_1} - \textbf{e}_{a_2}\right)}^TF_{\theta}(G)\left(\textbf{e}_{b_1} - \textbf{e}_{b_2}\right)\\
&\times \left(\cos\left(\frac{\Delta_\theta t}{2}\right) + \textbf{i}\frac{2-rn_2 - \theta}{\Delta_\theta}\sin\left(\frac{\Delta_\theta t}{2}\right)-1\right).
\end{align*}
Also, 
\begin{align*}
&\sum_{\theta \in \Phi_{a_1a_2} \setminus \{2r\}}\left|\frac{1}{2}\exp \left(\frac{-\textbf{i} \pi \theta}{2}\right){\left(\textbf{e}_{a_1} - \textbf{e}_{a_2}\right)}^TF_{\theta}(G)\left(\textbf{e}_{b_1} - \textbf{e}_{b_2}\right)\right|\\ 
&\leq \sum_{\theta \in \Phi_{a_1a_2}}\left|\frac{1}{2}\exp \left(\frac{-\textbf{i} \pi \theta}{2}\right){\left(\textbf{e}_{a_1} - \textbf{e}_{a_2}\right)}^TF_{\theta}(G)\left(\textbf{e}_{b_1} - \textbf{e}_{b_2}\right)\right|=1.
\end{align*}
The rest of the proof is similar to that of Theorem~\ref{Edge Corona First Theorem PGST}.
\end{proof}

\section{Infinite families of edge coronas admitting Pair-LPGST}\label{Infinite families}
In this section, we obtain several infinite families edge coronas admitting Pair-LPGST using Theorem~\ref{Edge Corona First Theorem PGST} and Theorem~\ref{Edge Corona Second Theorem PGST}. Further, we utilize  graph operations Cartesian product and tensor product to obtain a few more infinite families of edge coronas admitting Pair-LPGST. 

\noindent\textbf{Example 2.1.} Let $p$ be a prime such that $p\geq5$. Consider $G=\mathrm{Cay}(\mathbb{Z}_2 \oplus \mathbb{Z}_p,\{(1,0),(0,1),\ldots,(0,p-1)\})$. 
It can be verified that $G$ is a connected $r$-regular non-bipartite graph, where $r=p$. From Cao and Wan~\cite{XiwangCao2022} it follows that $G$ admits Pair-LPST at time $\frac{\pi}{2}$. Let $H=K_{n_2}$, where $n_2$ is an even integer such that $n_2+1 \geq p$ (For example $n_2=p-1$). Then it is clear that $(r-1)/(1+n_2)$ is not an integer. Now Theorem~\ref{Edge Corona First Theorem PGST} implies that the the edge corona $G \diamond H$ admits Pair-LPGST.   

\noindent\textbf{Example 2.2.} Let $m_1$ be an even integer. Consider $G=\overline{\cup_{m_1} K_2}$. It can be verified that $G$ is a connected $r$-regular non-bipartite graph, where $r=2m_1-2$. From Coutinho et al.~\cite{Coutinho2015} and Luo et al.~\cite{GaojunLuo2022}, it follows that $G$ admits Pair-LPST at time $\frac{\pi}{2}$. Let $H=K_{n_2}$, where $n_2$ is an even integer such that $n_2+1 \geq 2m_1-2$. Then it is clear that $(r-1)/(1+n_2)$ is not an integer. Now Theorem~\ref{Edge Corona First Theorem PGST} implies that the the edge corona $G \diamond H$ admits Pair-LPGST.   

\noindent\textbf{Example 2.3.} Let $n=4k$, where $k$ is an odd integer. Consider $G=\mathrm{Cay}(D_{2n}, \{a^k,a^{3k}\} \cup b \langle a \rangle)$. It can be verified that $G$ is a connected $r$-regular non-bipartite graph, where $r=4k+2$. From Cao and Feng~\cite{CaoDihedral2021} and Luo et al.~\cite{GaojunLuo2022}, it follows that $G$ admits Pair-LPST at time $\frac{\pi}{2}$. Let $H=K_{n_2}$, where $n_2$ is an even integer such that $n_2+1 \geq 4k+2$. Then it is clear that $(r-1)/(1+n_2)$ is not an integer. Now Theorem~\ref{Edge Corona First Theorem PGST} implies that the the edge corona $G \diamond H$ admits Pair-LPGST. 

\noindent\textbf{Example 2.4.} Let $n=2p$, where $p$ is an odd prime. Consider $G=\mathrm{Cay}(D_{2n}, \left(\langle a \rangle \setminus \{\textbf{1},a^p\}\right)\cup b \langle a \rangle)$. It can be verified that $G$ is a connected $r$-regular non-bipartite graph, where $r=4p-2$. From Cao and Feng~\cite{CaoDihedral2021} and Luo et al.~\cite{GaojunLuo2022}, it follows that $G$ admits Pair-LPST at time $\frac{\pi}{2}$. Let $H=K_{n_2}$, where $n_2$ is an even integer such that $n_2+1 \geq 4p-2$. Then it is clear that $(r-1)/(1+n_2)$ is not an integer. Now Theorem~\ref{Edge Corona First Theorem PGST} implies that the the edge corona $G \diamond H$ admits Pair-LPGST. 

\noindent\textbf{Example 2.5.} Let $n$ be an integer such that $n \geq 2$. Consider $G=\mathrm{Cay}(T_{4n}, T_{4n} \setminus \{\textbf{1},a^n\})$. It can be verified that $G$ is a connected $r$-regular non-bipartite graph, where $r=4n-2$. From Arezoomand et al.~\cite{ArezoomandDicyclic2022} and Luo et al.~\cite{GaojunLuo2022}, it follows that $G$ admits Pair-LPST at time $\frac{\pi}{2}$. Let $H=K_{n_2}$, where $n_2$ is an even integer such that $n_2+1 \geq 4n-2$. Then it is clear that $(r-1)/(1+n_2)$ is not an integer. Now Theorem~\ref{Edge Corona First Theorem PGST} implies that the the edge corona $G \diamond H$ admits Pair-LPGST. 

\noindent\textbf{Example 2.6.} Let $n$ be an odd integer such that $n \geq 3$. Consider $G=\mathrm{Cay}(SD_{8n}, \{a^n,a^{3n}\} \cup b \langle a \rangle)$. It can be verified that $G$ is a connected $r$-regular non-bipartite graph, where $r=4n+2$. From Luo et al.~\cite{LuoSemidihedral2022} and Luo et al.~\cite{GaojunLuo2022}, it follows that $G$ admits Pair-LPST at time $\frac{\pi}{2}$. Let $H=K_{n_2}$, where $n_2$ is an even integer such that $n_2+1 \geq 4n+2$. Then it is clear that $(r-1)/(1+n_2)$ is not an integer. Now Theorem~\ref{Edge Corona First Theorem PGST} implies that the the edge corona $G \diamond H$ admits Pair-LPGST.

\noindent\textbf{Example 2.7.} Let $r$ be an odd integer such that $r \geq 3$. Consider $G=Q_r$ and $H=K_{4r}$. Clearly, $G$ is a connected $r$-regular bipartite graph. It follows from Christandl et al.~\cite{Christandl2004} and Luo et al.~\cite{GaojunLuo2022} that $G$ admits Pair-LPST at time $\frac{\pi}{2}$. Also, $n_2=4r$, $r(n_2-2) \equiv 2~(\mathrm{mod}~4)$ and $(r-1)/(1+n_2)$ is not an integer. Now Theorem~\ref{Edge Corona Second Theorem PGST} implies that the the edge corona $G \diamond H$ admits Pair-LPGST. 

Let $G_1$ and $G_2$ be two graphs. The \textit{Cartesian product} of $G_1$ and $G_2$ is denoted by $G_1 \square G_2$ and is defined as the graph having vertex set $V(G_1) \times V(G_2)$ such that two vertices $(g_1,g_2)$ and $(h_1,h_2)$ are adjacent in $G_1 \square G_2$ if $g_1=h_1$ and $g_2$ is adjacent to $h_2$ in $G_2$; or $g_1$ is adjacent to $h_1$ in $G_1$ and $g_2=h_2$. The \textit{tensor product} of $G_1$ and $G_2$ is denoted by $G_1 \times G_2$ and is defined as the graph having vertex set $V(G_1) \times V(G_2)$ such that two vertices $(g_1,g_2)$ and $(h_1,h_2)$ are adjacent in $G_1 \times G_2$ if $g_1$ is adjacent to $h_1$ in $G_1$ and $g_2$ is adjacent to $h_2$ in $G_2$. 

In the following, we state some useful results.
\begin{theorem}\emph{\cite{CayleyGraphProduct2016}}\label{Cayley Graph Cartesian Product}
Let $\mathrm{Cay}(\Gamma_1,S_1)$ and $\mathrm{Cay}(\Gamma_2,S_2)$ be two Cayley graphs. Then $\mathrm{Cay}(\Gamma_1,S_1) \square \mathrm{Cay}(\Gamma_2,S_2)$ is also a Cayley graph given by $\mathrm{Cay}(\Gamma_1,S_1) \square \mathrm{Cay}(\Gamma_2,S_2) = \mathrm{Cay}(\Gamma_1 \times \Gamma_2, (S_1 \times \textbf{1}) \cup ({\textbf{1} \times S_2}))$.
\end{theorem}

\begin{theorem}\emph{\cite{CayleyGraphProduct2016}}\label{Cayley Graph Tensor Product}
Let $\mathrm{Cay}(\Gamma_1,S_1)$ and $\mathrm{Cay}(\Gamma_2,S_2)$ be two Cayley graphs. Then $\mathrm{Cay}(\Gamma_1,S_1) \times \mathrm{Cay}(\Gamma_2,S_2)$ is also a Cayley graph given by $\mathrm{Cay}(\Gamma_1,S_1) \times \mathrm{Cay}(\Gamma_2,S_2) = \mathrm{Cay}(\Gamma_1 \times \Gamma_2, S_1 \times S_2)$. 
\end{theorem}

\begin{theorem}\emph{\cite{Canul2010}}\label{PST and Cartesian product}
Let $G_1$ and $G_2$ be two graphs such that both admit vertex PST at time $t$. Then the graph $G_1 \square G_2$ also admits vertex PST at time $t$.
\end{theorem}

\begin{theorem}\emph{\cite{GraphProduct2011}}\label{PST and tensor product}
Let $G_1$ be a graph that admit vertex PST at time $t$ and $t\lambda = m_1 \pi$, for some $m_1 \in \mathbb{Z}$ and for all adjacency eigenvalues $\lambda$ of $G_1$. Also, let $G_2$ be a circulant graph such that all the adjacency eigenvalues are odd. Then the graph $G_1 \times G_2$ admits vertex PST at time $t$.
\end{theorem}
    
Now we use Cartesian product and tensor product to obtain more graphs admitting Pair-LPGST in the next two examples. 

\noindent\textbf{Example 2.8.} Let $\mathrm{Cay}(\Gamma_1,S_1)$ and $\mathrm{Cay}(\Gamma_2,S_2)$ be two connected Cayley graphs, where $\mathrm{Cay}(\Gamma_1,S_1)$ is non-bipartite. Also, let both of $\mathrm{Cay}(\Gamma_1,S_1)$ and $\mathrm{Cay}(\Gamma_2,S_2)$ admit vertex PST at time $\frac{\pi}{2}$. Then the graph $G$ is $r$-regular connected non-bipartite, where $G = \mathrm{Cay}(\Gamma_1,S_1) \square \mathrm{Cay}(\Gamma_2,S_2)$ and $r = |S_1| + |S_2|$. Now Theorem~\ref{Cayley Graph Cartesian Product} and Theorem 3.1 of Luo et al.~\cite{GaojunLuo2022} together implies that $G$ admits Pair-LPST at time $\frac{\pi}{2}$. Consider $H=K_{n_2}$, where $n_2$ is an even integer such that $n_2+1 \geq r$ (For example $n_2=2r$). It is clear that $(r-1)/(1+n_2)$ is not an integer. Therefore, Theorem~\ref{Edge Corona First Theorem PGST} implies that $G \diamond H$ admits Pair-LPGST.
      
\noindent\textbf{Example 2.9.} Let $\mathrm{Cay}(\Gamma,S)$ be a connected non-bipartite Cayley graph such that all the adjacency eigenvalues are even. Let $\mathrm{Cay}(\Gamma,S)$ admit vertex PST at time $\frac{\pi}{2}$. Also, let $m$ be an even integer such that $m \geq 4$. Then $G$ is $r$-regular connected non-bipartite, where $G = \mathrm{Cay}(\Gamma,S) \times K_m$ and $r = |S|(m-1)$. Now Theorem~\ref{Cayley Graph Tensor Product} and Theorem 3.1 of Luo et al.~\cite{GaojunLuo2022} together implies that $G$ admits Pair-LPST at time $\frac{\pi}{2}$. If $H=K_{n_2}$ such that $n_2$ is an even integer and $n_2+1 \geq r$, then $(r-1)/(1+n_2)$ is not an integer. Therefore, Theorem~\ref{Edge Corona First Theorem PGST} implies that $G \diamond H$ admits Pair-LPGST. 
      
Note that the method used in Example 2.8 and Example 2.9, and utilizing the Cayley graphs appearing in the beginning of this section, one can construct several new infinite families of edge coronas admitting Pair-LPGST.

From Theorem 6.1 of Godsil et al.~\cite{RealStateTransfer2025}, it follows that the complete graph does not admit Pair-LPST. Therefore, Theorem~\ref{Edge Corona First Theorem PGST}, Theorem~\ref{Edge Corona Second Theorem PGST}, Theorem~\ref{Graphs with clusters Pair-LPST}, Theorem~\ref{non-existence of Pair-LPST First Theorem} and Theorem~\ref{non-existence of Pair-LPST 2nd Lemma} implies that the edge coronas appearing in Example 1 to Example 9 provides infinitely many graphs admitting Pair-LPGST that fails to admit Pair-LPST.  
      
\section{Pair-LPGST on edge coronas of complete graphs}\label{Edge corona of the complete graph}
In this section, we investigate the existence of Pair-LPGST on the edge corona $K_{n_1}\diamond K_{n_2}$, where $n_1$ and $n_2$ are two positive integer such that $n_1 \geq 4$. The eigenvalues of the graph $K_{n_1}$ are given by $\theta_0=0$ and $\theta_1=n_1$. Further, the eigenprojectors corresponding to $\theta_0$ and $\theta_1$  are given by 
$$F_{\theta_0}(K_{n_1}) = \frac{1}{n_1}J_{n_1 \times n_1}~~\mathrm{and}~~F_{\theta_1}(K_{n_1}) = I_{n_1 \times n_1} - \frac{1}{n_1}J_{n_1 \times n_1},$$
respectively. In the following lemma, we prove that a pair $((a_1,0),(a_2,0))$ of $K_{n_1}\diamond K_{n_2}$, where $a_1$ and $a_2$ are vertices of $K_{n_1}$ cannot be involved in Pair-LPGST.

\begin{lema}
Let $a_1,a_2,b_1$ and $b_2$ be vertices of $K_{n_1}$. Then the edge corona $K_{n_1}\diamond K_{n_2}$ cannot admits Pair-LPGST from the pair $((a_1,0),(a_2,0))$ to $((b_1,0),(b_2,0))$.
\end{lema}
\begin{proof}
If possible, let $K_{n_1}\diamond K_{n_2}$ admit Pair-LPGST from $((a_1,0),(a_2,0))$ to $((b_1,0),(b_2,0))$. Then it follows from Bommel~\cite{Bommel} that $((a_1,0),(a_2,0))$ and $((b_1,0),(b_2,0))$ are strongly cospectral. From Theorem~\ref{Non-bipartite Spectrum}, for each eigenvalue $\mu$ of $K_{n_2}$, we obtain that
$$F_{\mu+2}(G \diamond \mathcal{H})(\textbf{e}_{(a_1,0)} - \textbf{e}_{(a_2,0)}) = \textbf{0}_{n_1+mn_2} = F_{\mu+2}(G \diamond \mathcal{H})(\textbf{e}_{(b_1,0)} - \textbf{e}_{(b_2,0)}).$$
Further, we obtain that
$$F_{{\theta_0}^{\pm}}(G \diamond \mathcal{H})(\textbf{e}_{(a_1,0)} - \textbf{e}_{(a_2,0)}) = \textbf{0}_{n_1+mn_2} = F_{{\theta_0}^{\pm}}(G \diamond \mathcal{H})(\textbf{e}_{(b_1,0)} - \textbf{e}_{(b_2,0)}).$$
This implies that $\Phi_{(a_1,0)(a_2,0)} \subseteq \{{\theta_1}^{+},{\theta_1}^{-}\}$.
Since $((a_1,0),(a_2,0))$ and $((b_1,0),(b_2,0))$ are strongly cospectral, Lemma~\ref{fixed state} implies that $\Phi_{(a_1,0)(a_2,0)} = \{{\theta_1}^{+},{\theta_1}^{-}\}$. Then by Lemma~\ref{Laplacian perfect pair state transfer lemma non-bipartite}, we have either ${\theta_1}^{+},{\theta_1}^{-} \in \Phi^{+}_{(a_1,0)(a_2,0),(b_1,0)(b_2,0)}$ or ${\theta_1}^{+},{\theta_1}^{-} \in \Phi^{-}_{(a_1,0)(a_2,0),(b_1,0)(b_2,0)}$, that is, either $\Phi^{-}_{(a_1,0)(a_2,0),(b_1,0)(b_2,0)}$ is empty or $\Phi^{+}_{(a_1,0)(a_2,0),(b_1,0)(b_2,0)}$ is empty, which contradicts Lemma~\ref{fixed state}. Therefore, $K_{n_1}\diamond K_{n_2}$ cannot admits Pair-LPGST from $((a_1,0),(a_2,0))$ to $((b_1,0),(b_2,0))$.
\end{proof}
Since strong cospectrality is a necessary condition for the existence of Pair-LPGST, it follows from Lemma~\ref{fixed state} that the graph $K_{n_2}$ does not admit Pair-LPGST. Therefore, the next lemma follows directly from Theorem~\ref{Graphs with clusters Pair-LPGST}.
\begin{lema}\emph{\cite{MonterdePal2026}}
Let $e_s$ be an edge of $K_{n_1}$ and $a_1,a_2,b_1$, and $b_2$ be vertices of $K_{n_2}$. Then the edge corona $K_{n_1}\diamond K_{n_2}$ cannot admits Pair-LPGST from the pair $((a_1,e_s),(a_2,e_s))$ to $((b_1,e_s),(b_2,e_s))$.
\end{lema}

\begin{lema}
Let $a_1$ and $b_1$ be two distinct vertices of $K_{n_1}$. Also, let $e_s$ be an edge of $K_{n_1}$ and $a_2,b_2$ be vertices of $K_{n_2}$. Then the edge corona $K_{n_1}\diamond K_{n_2}$ cannot admits Pair-LPGST from the pair $((a_1,0),(a_2,e_s))$ to $((b_1,0),(b_2,e_s))$.
\end{lema}
\begin{proof}
If possible, let $K_{n_1}\diamond K_{n_2}$ admit Pair-LPGST from $((a_1,0),(a_2,e_s))$ to $((b_1,0),(b_2,e_s))$. Then $((a_1,0),(a_2,0))$ and $((b_1,0),(b_2,0))$ are strongly cospectral. This implies
$$F_{\theta^+_1}(G \diamond \mathcal{H})(\textbf{e}_{(a_1,0)} - \textbf{e}_{(a_2,e_s)}) = \pm F_{\theta^+_1}(G \diamond \mathcal{H})(\textbf{e}_{(b_1,0)} - \textbf{e}_{(b_2,e_s)}).$$
Using Theorem~\ref{Non-bipartite Spectrum}, we obtain that
\begin{align*}
&F_{\theta^+_1}(G \diamond \mathcal{H})(\textbf{e}_{(a_1,0)}-\textbf{e}_{(a_2,e_s)})\\ &= \frac{(2-\theta^+_1)^2}{(2-\theta^+_1)^2+(2r-\theta_1)n_2}\begin{pmatrix}
     F_{\theta_1}(K_{n_1})\textbf{e}_{a_1}-\frac{1}{2-\theta^+_1}F_{\theta_1}(K_{n_1})M_{K_{n_1}}\textbf{e}_s\\
         \frac{1}{2-\theta^+_1}M^T_{K_{n_1}}F_{\theta_1}(K_{n_1})\textbf{e}_{a_1} \otimes \textbf{1}_{n_2} - \frac{1}{{(2-\theta^+_1)}^2}M^T_{K_{n_1}}F_{\theta_1}(K_{n_1})M_{K_{n_1}}\textbf{e}_s\otimes \textbf{1}_{n_2}
         \end{pmatrix}
\end{align*}
and 
\begin{align*}
&F_{\theta^+_1}(G \diamond \mathcal{H})(\textbf{e}_{(b_1,0)}-\textbf{e}_{(b_2,e_s)})\\ &= \frac{(2-\theta^+_1)^2}{(2-\theta^+_1)^2+(2r-\theta_1)n_2}\begin{pmatrix}
     F_{\theta_1}(K_{n_1})\textbf{e}_{b_1}-\frac{1}{2-\theta^+_1}F_{\theta_1}(K_{n_1})M_{K_{n_1}}\textbf{e}_s\\
         \frac{1}{2-\theta^+_1}M^T_{K_{n_1}}F_{\theta_1}(K_{n_1})\textbf{e}_{b_1} \otimes \textbf{1}_{n_2} - \frac{1}{{(2-\theta^+_1)}^2}M^T_{K_{n_1}}F_{\theta_1}(K_{n_1})M_{K_{n_1}}\textbf{e}_s\otimes \textbf{1}_{n_2}
         \end{pmatrix}.
\end{align*}
The following two cases arise.

\noindent \textbf{Case 1.} $F_{\theta^+_1}(G \diamond \mathcal{H})(\textbf{e}_{(a_1,0)} - \textbf{e}_{(a_2,e_s)}) = F_{\theta^+_1}(G \diamond \mathcal{H})(\textbf{e}_{(b_1,0)} - \textbf{e}_{(b_2,e_s)})$. This implies  
$$F_{\theta_1}(K_{n_1})\textbf{e}_{a_1} = F_{\theta_1}(K_{n_1})\textbf{e}_{b_1}.$$
Recall that $F_{\theta_1}(K_{n_1}) = I_{n_1 \times n_1} - \frac{1}{n_1}J_{n_1 \times n_1}$. Therefore, we have $\textbf{e}_{a_1} - \frac{1}{n_1}\textbf{1}_{n_1} = \textbf{e}_{b_1} - \frac{1}{n_1}\textbf{1}_{n_1}$. This implies that $a_1=b_1$, a contradiction. 

\noindent \textbf{Case 2.} $F_{\theta^+_1}(G \diamond \mathcal{H})(\textbf{e}_{(a_1,0)} - \textbf{e}_{(a_2,e_s)}) = -F_{\theta^+_1}(G \diamond \mathcal{H})(\textbf{e}_{(b_1,0)} - \textbf{e}_{(b_2,e_s)})$. This implies that 
\begin{equation*}
F_{\theta_1}(K_{n_1})\textbf{e}_{a_1} + F_{\theta_1}(K_{n_1})\textbf{e}_{b_1} = \frac{2}{2-\theta^+_1}F_{\theta_1}(K_{n_1})M_{K_{n_1}}\textbf{e}_s.
\end{equation*}
This further implies that 
\begin{equation}\label{d_7}
\textbf{e}_{a_1} + \textbf{e}_{b_1} - \frac{2}{n_1}\textbf{1}_{n_1}= \frac{2}{2-\theta^+_1}I_{n_1 \times n_1}M_{K_{n_1}}\textbf{e}_s - \frac{4}{n_1(2-\theta^+_1)}\textbf{1}_{n_1}.
\end{equation}
If $v_i$ and $v_j$ are the end vertices of the edge $e_s$, then it can be verified that $I_{n_1 \times n_1}M_{K_{n_1}}\textbf{e}_s = \textbf{e}_{v_i} + \textbf{e}_{v_j}$. Thus Equation~(\ref{d_7}) reduces to
\begin{equation}\label{d_8}
\textbf{e}_{a_1} + \textbf{e}_{b_1} - \frac{2}{n_1}\textbf{1}_{n_1}= \frac{2}{2-\theta^+_1}\left(\textbf{e}_{v_i} + \textbf{e}_{v_j} \right) - \frac{4}{n_1(2-\theta^+_1)}\textbf{1}_{n_1}.
\end{equation}
Now we consider the following sub-cases.

\noindent \textbf{Sub-Case 1.} If both of $a_1$ and $b_1$ are end vertices of $e_s$, then we obtain that
\begin{equation}
\textbf{e}_{a_1} + \textbf{e}_{b_1} - \frac{2}{n_1}\textbf{1}_{n_1}= \frac{2}{2-\theta^+_1}\left(\textbf{e}_{a_1} + \textbf{e}_{b_1} \right) - \frac{4}{n_1(2-\theta^+_1)}\textbf{1}_{n_1}.
\end{equation}
Pre-multiplying both sides of the previous equation by $\textbf{e}^T_{a_1}$, we have 
$$1-\frac{2}{n_1} = \frac{2}{2-\theta^+_1}\left(1-\frac{2}{n_1}\right).$$
This implies that $\theta^+_1 = 0$, a contradiction.

\noindent \textbf{Sub-Case 2.} If exactly one of $a_1$ and $b_1$ is an end vertex of $e_s$, then proceeding as in the previous sub-case, we arrive at a contradiction.

\noindent \textbf{Sub-Case 3.} If neither $a_1$ nor $b_1$ is an end vertex of $e_s$, then pre-multiplying both sides of the Equation~(\ref{d_8}) by $\textbf{e}^T_{a_1}$, we obtain that 
$$1-\frac{2}{n_1} = \frac{4}{n_1(\theta^+_1-2)}.$$
Note that $\frac{4}{n_1(\theta^+_1-2)} < \frac{2}{n_1}$. Therefore, the preceding equation implies that $n_1 < 4$, a contradiction to the fact that  $n_1 \geq 4$. Therefore, the edge corona $K_{n_1}\diamond K_{n_2}$ cannot admits Pair-LPGST from the pair $((a_1,0),(a_2,e_s))$ to $((b_1,0),(b_2,e_s))$. This completes the proof.
\end{proof}

\subsection*{Acknowledgements} The first author acknowledges the
support provided by the Prime Minister’s Research Fellowship (PMRF), PMRF-ID:
1903283, Government of India.

\end{document}